\documentstyle[12pt]{article}

\makeatletter
\def\@eqnnum{\hbox to .01pt{}\rlap{\rm \hskip -\displaywidth(\theequation)}}
\makeatother

\makeatletter

\font\fivmsam=msam5
\font\sixmsam=msam6
\font\sevenmsam=msam7
\font\egtmsam=msam8
\font\ninmsam=msam9
\font\tenmsam=msam10
\font\elvmsam=msam10 \@halfmag
\font\twlmsam=msam10 \@magscale1
\font\frtnmsam=msam10 \@magscale2
\font\svtnmsam=msam10 \@magscale3
\font\twtymsam=msam10 \@magscale4
\font\twfvmsam=msam10 \@magscale5
\font\fivmsbm=msbm5
\font\sixmsbm=msbm6
\font\sevenmsbm=msbm7
\font\egtmsbm=msbm8
\font\ninmsbm=msbm9
\font\tenmsbm=msbm10
\font\elvmsbm=msbm10 \@halfmag
\font\twlmsbm=msbm10 \@magscale1
\font\frtnmsbm=msbm10 \@magscale2
\font\svtnmsbm=msbm10 \@magscale3
\font\twtymsbm=msbm10 \@magscale4
\font\twfvmsbm=msbm10 \@magscale5
\font\fiveul=eufm5
\font\sixeul=eufm6
\font\seveneul=eufm7
\font\egteul=eufm8
\font\nineul=eufm9
\font\teneul=eufm10
\font\elveul=eufm10 \@halfmag
\font\twleul=eufm10 \@magscale1
\font\frtneul=eufm10 \@magscale2
\font\svtneul=eufm10 \@magscale3
\font\twtyeul=eufm10 \@magscale4
\font\twfveul=eufm10 \@magscale5
\font\egtcmssi=cmssi8
\font\nincmssi=cmssi9
\font\tencmssi=cmssi10
 \@halfmag

\newfam\msamfam
\newfam\msbmfam
\newfam\eulfam
\newfam\sfifam

\def\@vpt{\def\msam{\fam\msamfam\fivmsam}%
  \textfont\msamfam\fivmsam
  \scriptfont\msamfam\fivmsam
  \scriptscriptfont\msamfam\fivmsam
  \def\msbm{\fam\msbmfam\fivmsbm}%
  \textfont\msbmfam\fivmsbm
  \scriptfont\msbmfam\fivmsbm
  \scriptscriptfont\msbmfam\fivmsbm
  \def\eul{\fam\eulfam\fiveul}%
  \textfont\eulfam\fiveul
  \scriptfont\eulfam\fiveul
  \scriptscriptfont\eulfam\fiveul}%
\def\@vipt{\def\msam{\fam\msamfam\sixmsam}%
  \textfont\msamfam\sixmsam
  \scriptfont\msamfam\sixmsam
  \scriptscriptfont\msamfam\sixmsam
  \def\msbm{\fam\msbmfam\sixmsbm}%
  \textfont\msbmfam\sixmsbm
  \scriptfont\msbmfam\sixmsbm
  \scriptscriptfont\msbmfam\sixmsbm
  \def\eul{\fam\eulfam\sixeul}%
  \textfont\eulfam\sixeul
  \scriptfont\eulfam\sixeul
  \scriptscriptfont\eulfam\sixeul}%
\def\@viipt{\def\msam{\fam\msamfam\sevenmsam}%
  \textfont\msamfam\sevenmsam
  \scriptfont\msamfam\fivmsam
  \scriptscriptfont\msamfam\fivmsam
  \def\msbm{\fam\msbmfam\sevenmsbm}%
  \textfont\msbmfam\sevenmsbm
  \scriptfont\msbmfam\fivmsbm
  \scriptscriptfont\msbmfam\fivmsbm
  \def\eul{\fam\eulfam\seveneul}%
  \textfont\eulfam\seveneul
  \scriptfont\eulfam\fiveul
  \scriptscriptfont\eulfam\fiveul}%
\def\@viiipt{\def\msam{\fam\msamfam\egtmsam}%
  \textfont\msamfam\egtmsam
  \scriptfont\msamfam\sixmsam
  \scriptscriptfont\msamfam\fivmsam
  \def\msbm{\fam\msbmfam\egtmsbm}%
  \textfont\msbmfam\egtmsbm
  \scriptfont\msbmfam\sixmsbm
  \scriptscriptfont\msbmfam\fivmsbm
  \def\eul{\fam\eulfam\egteul}%
  \textfont\eulfam\egteul
  \scriptfont\eulfam\sixeul
  \scriptscriptfont\eulfam\fiveul
  \def\sfi{\fam\sfifam\egtcmssi}%
  \textfont\sfifam\egtcmssi
  \scriptfont\sfifam\egtcmssi
  \scriptscriptfont\sfifam\egtcmssi}%
\def\@ixpt{\def\msam{\fam\msamfam\ninmsam}%
  \textfont\msamfam\ninmsam
  \scriptfont\msamfam\sixmsam
  \scriptscriptfont\msamfam\fivmsam
  \def\msbm{\fam\msbmfam\ninmsbm}%
  \textfont\msbmfam\ninmsbm
  \scriptfont\msbmfam\sixmsbm
  \scriptscriptfont\msbmfam\fivmsbm
  \def\eul{\fam\eulfam\nineul}%
  \textfont\eulfam\nineul
  \scriptfont\eulfam\sixeul
  \scriptscriptfont\eulfam\fiveul
  \def\sfi{\fam\sfifam\nincmssi}%
  \textfont\sfifam\nincmssi
  \scriptfont\sfifam\egtcmssi
  \scriptscriptfont\sfifam\egtcmssi}%
\def\@xpt{\def\msam{\fam\msamfam\tenmsam}%
  \textfont\msamfam\tenmsam
  \scriptfont\msamfam\sevenmsam
  \scriptscriptfont\msamfam\fivmsam
  \def\msbm{\fam\msbmfam\tenmsbm}%
  \textfont\msbmfam\tenmsbm
  \scriptfont\msbmfam\sevenmsbm
  \scriptscriptfont\msbmfam\fivmsbm
  \def\eul{\fam\eulfam\teneul}%
  \textfont\eulfam\teneul
  \scriptfont\eulfam\seveneul
  \scriptscriptfont\eulfam\fiveul
  \def\sfi{\fam\sfifam\tencmssi}%
  \textfont\sfifam\tencmssi
  \scriptfont\sfifam\egtcmssi
  \scriptscriptfont\sfifam\egtcmssi}%
\def\@xipt{\def\msam{\fam\msamfam\elvmsam}%
  \textfont\msamfam\elvmsam
  \scriptfont\msamfam\egtmsam
  \scriptscriptfont\msamfam\sixmsam
  \def\msbm{\fam\msbmfam\elvmsbm}%
  \textfont\msbmfam\elvmsbm
  \scriptfont\msbmfam\egtmsbm
  \scriptscriptfont\msbmfam\sixmsbm
  \def\eul{\fam\eulfam\elveul}%
  \textfont\eulfam\elveul
  \scriptfont\eulfam\egteul
  \scriptscriptfont\eulfam\sixeul}%
\def\@xiipt{\def\msam{\fam\msamfam\twlmsam}%
  \textfont\msamfam\twlmsam
  \scriptfont\msamfam\egtmsam
  \scriptscriptfont\msamfam\sixmsam
  \def\msbm{\fam\msbmfam\twlmsbm}%
  \textfont\msbmfam\twlmsbm
  \scriptfont\msbmfam\egtmsbm
  \scriptscriptfont\msbmfam\sixmsbm
  \def\eul{\fam\eulfam\twleul}%
  \textfont\eulfam\twleul
  \scriptfont\eulfam\egteul
  \scriptscriptfont\eulfam\sixeul}%
\def\@xivpt{\def\msam{\fam\msamfam\frtnmsam}%
  \textfont\msamfam\frtnmsam
  \scriptfont\msamfam\tenmsam
  \scriptscriptfont\msamfam\sevenmsam
  \def\msbm{\fam\msbmfam\frtnmsbm}%
  \textfont\msbmfam\frtnmsbm
  \scriptfont\msbmfam\tenmsbm
  \scriptscriptfont\msbmfam\sevenmsbm
  \def\eul{\fam\eulfam\frtneul}%
  \textfont\eulfam\frtneul
  \scriptfont\eulfam\teneul
  \scriptscriptfont\eulfam\seveneul}%
\def\@xviipt{\def\msam{\fam\msamfam\svtnmsam}%
  \textfont\msamfam\svtnmsam
  \scriptfont\msamfam\twlmsam
  \scriptscriptfont\msamfam\tenmsam
  \def\msbm{\fam\msbmfam\svtnmsbm}%
  \textfont\msbmfam\svtnmsbm
  \scriptfont\msbmfam\twlmsbm
  \scriptscriptfont\msbmfam\tenmsbm
  \def\eul{\fam\eulfam\svtneul}%
  \textfont\eulfam\svtneul
  \scriptfont\eulfam\twleul
  \scriptscriptfont\eulfam\teneul}%
\def\@xxpt{\def\msam{\fam\msamfam\twtymsam}%
  \textfont\msamfam\twtymsam
  \scriptfont\msamfam\frtnmsam
  \scriptscriptfont\msamfam\twlmsam
  \def\msbm{\fam\msbmfam\twtymsbm}%
  \textfont\msbmfam\twtymsbm
  \scriptfont\msbmfam\frtnmsbm
  \scriptscriptfont\msbmfam\twlmsbm
  \def\eul{\fam\eulfam\twtyeul}%
  \textfont\eulfam\twtyeul
  \scriptfont\eulfam\frtneul
  \scriptscriptfont\eulfam\twleul}%
\def\@xxvpt{\def\msam{\fam\msamfam\twfvmsam}%
  \textfont\msamfam\twfvmsam
  \scriptfont\msamfam\twtymsam
  \scriptscriptfont\msamfam\svtnmsam
  \def\msbm{\fam\msbmfam\twfvmsbm}%
  \textfont\msbmfam\twfvmsbm
  \scriptfont\msbmfam\twtymsbm
  \scriptscriptfont\msbmfam\svtnmsbm
  \def\eul{\fam\eulfam\twfveul}%
  \textfont\eulfam\twfveul
  \scriptfont\eulfam\twtyeul
  \scriptscriptfont\eulfam\svtneul}%

\def\hexnumber@#1{\ifnum#1<10 \number#1\else
 \ifnum#1=10 A\else\ifnum#1=11 B\else\ifnum#1=12 C\else
 \ifnum#1=13 D\else\ifnum#1=14 E\else\ifnum#1=15 F\fi\fi\fi\fi\fi\fi\fi}%

\def\msa@{\hexnumber@\msamfam}%
\def\msb@{\hexnumber@\msbmfam}%
\mathchardef\boxdot="2\msa@00
\mathchardef\boxplus="2\msa@01
\mathchardef\boxtimes="2\msa@02
\mathchardef\square="0\msa@03
\mathchardef\blacksquare="0\msa@04
\mathchardef\centerdot="2\msa@05
\mathchardef\lozenge="0\msa@06
\mathchardef\blacklozenge="0\msa@07
\mathchardef\circlearrowright="3\msa@08
\mathchardef\circlearrowleft="3\msa@09
\mathchardef\rightleftharpoons="3\msa@0A
\mathchardef\leftrightharpoons="3\msa@0B
\mathchardef\boxminus="2\msa@0C
\mathchardef\Vdash="3\msa@0D
\mathchardef\Vvdash="3\msa@0E
\mathchardef\vDash="3\msa@0F
\mathchardef\twoheadrightarrow="3\msa@10
\mathchardef\twoheadleftarrow="3\msa@11
\mathchardef\leftleftarrows="3\msa@12
\mathchardef\rightrightarrows="3\msa@13
\mathchardef\upuparrows="3\msa@14
\mathchardef\downdownarrows="3\msa@15
\mathchardef\upharpoonright="3\msa@16
\mathchardef\downharpoonright="3\msa@17
\mathchardef\upharpoonleft="3\msa@18
\mathchardef\downharpoonleft="3\msa@19
\mathchardef\rightarrowtail="3\msa@1A
\mathchardef\leftarrowtail="3\msa@1B
\mathchardef\leftrightarrows="3\msa@1C
\mathchardef\rightleftarrows="3\msa@1D
\mathchardef\Lsh="3\msa@1E
\mathchardef\Rsh="3\msa@1F
\mathchardef\rightsquigarrow="3\msa@20
\mathchardef\leftrightsquigarrow="3\msa@21
\mathchardef\looparrowleft="3\msa@22
\mathchardef\looparrowright="3\msa@23
\mathchardef\circeq="3\msa@24
\mathchardef\succsim="3\msa@25
\mathchardef\gtrsim="3\msa@26
\mathchardef\gtrapprox="3\msa@27
\mathchardef\multimap="3\msa@28
\mathchardef\therefore="3\msa@29
\mathchardef\because="3\msa@2A
\mathchardef\doteqdot="3\msa@2B

\mathchardef\triangleq="3\msa@2C
\mathchardef\precsim="3\msa@2D
\mathchardef\lesssim="3\msa@2E
\mathchardef\lessapprox="3\msa@2F
\mathchardef\eqslantless="3\msa@30
\mathchardef\eqslantgtr="3\msa@31
\mathchardef\curlyeqprec="3\msa@32
\mathchardef\curlyeqsucc="3\msa@33
\mathchardef\preccurlyeq="3\msa@34
\mathchardef\leqq="3\msa@35
\mathchardef\leqslant="3\msa@36
\mathchardef\lessgtr="3\msa@37
\mathchardef\backprime="0\msa@38
\mathchardef\risingdotseq="3\msa@3A
\mathchardef\fallingdotseq="3\msa@3B
\mathchardef\succcurlyeq="3\msa@3C
\mathchardef\geqq="3\msa@3D
\mathchardef\geqslant="3\msa@3E
\mathchardef\gtrless="3\msa@3F
\mathchardef\sqsubset="3\msa@40
\mathchardef\sqsupset="3\msa@41
\mathchardef\vartriangleright="3\msa@42
\mathchardef\vartriangleleft="3\msa@43
\mathchardef\trianglerighteq="3\msa@44
\mathchardef\trianglelefteq="3\msa@45
\mathchardef\bigstar="0\msa@46
\mathchardef\between="3\msa@47
\mathchardef\blacktriangledown="0\msa@48
\mathchardef\blacktriangleright="3\msa@49
\mathchardef\blacktriangleleft="3\msa@4A
\mathchardef\vartriangle="3\msa@4D
\mathchardef\blacktriangle="0\msa@4E
\mathchardef\triangledown="0\msa@4F
\mathchardef\eqcirc="3\msa@50
\mathchardef\lesseqgtr="3\msa@51
\mathchardef\gtreqless="3\msa@52
\mathchardef\lesseqqgtr="3\msa@53
\mathchardef\gtreqqless="3\msa@54
\mathchardef\Rrightarrow="3\msa@56
\mathchardef\Lleftarrow="3\msa@57
\mathchardef\veebar="2\msa@59
\mathchardef\barwedge="2\msa@5A
\mathchardef\doublebarwedge="2\msa@5B
\mathchardef\angle="0\msa@5C
\mathchardef\measuredangle="0\msa@5D
\mathchardef\sphericalangle="0\msa@5E
\mathchardef\varpropto="3\msa@5F
\mathchardef\smallsmile="3\msa@60
\mathchardef\smallfrown="3\msa@61
\mathchardef\Subset="3\msa@62
\mathchardef\Supset="3\msa@63
\mathchardef\Cup="2\msa@64

\mathchardef\Cap="2\msa@65

\mathchardef\curlywedge="2\msa@66
\mathchardef\curlyvee="2\msa@67
\mathchardef\leftthreetimes="2\msa@68
\mathchardef\rightthreetimes="2\msa@69
\mathchardef\subseteqq="3\msa@6A
\mathchardef\supseteqq="3\msa@6B
\mathchardef\bumpeq="3\msa@6C
\mathchardef\Bumpeq="3\msa@6D
\mathchardef\lll="3\msa@6E

\mathchardef\ggg="3\msa@6F

\mathchardef\circledS="0\msa@73
\mathchardef\pitchfork="3\msa@74
\mathchardef\dotplus="2\msa@75
\mathchardef\backsim="3\msa@76
\mathchardef\backsimeq="3\msa@77
\mathchardef\complement="0\msa@7B
\mathchardef\intercal="2\msa@7C
\mathchardef\circledcirc="2\msa@7D
\mathchardef\circledast="2\msa@7E
\mathchardef\circleddash="2\msa@7F
\def\ulcorner{\delimiter"4\msa@70\msa@70 }%
\def\urcorner{\delimiter"5\msa@71\msa@71 }%
\def\llcorner{\delimiter"4\msa@78\msa@78 }%
\def\lrcorner{\delimiter"5\msa@79\msa@79 }%
\def\yen{\mathhexbox\msa@55 }%
\def\checkmark{\mathhexbox\msa@58 }%
\def\circledR{\mathhexbox\msa@72 }%
\def\maltese{\mathhexbox\msa@7A }%
\mathchardef\lvertneqq="3\msb@00
\mathchardef\gvertneqq="3\msb@01
\mathchardef\nleq="3\msb@02
\mathchardef\ngeq="3\msb@03
\mathchardef\nless="3\msb@04
\mathchardef\ngtr="3\msb@05
\mathchardef\nprec="3\msb@06
\mathchardef\nsucc="3\msb@07
\mathchardef\lneqq="3\msb@08
\mathchardef\gneqq="3\msb@09
\mathchardef\nleqslant="3\msb@0A
\mathchardef\ngeqslant="3\msb@0B
\mathchardef\lneq="3\msb@0C
\mathchardef\gneq="3\msb@0D
\mathchardef\npreceq="3\msb@0E
\mathchardef\nsucceq="3\msb@0F
\mathchardef\precnsim="3\msb@10
\mathchardef\succnsim="3\msb@11
\mathchardef\lnsim="3\msb@12
\mathchardef\gnsim="3\msb@13
\mathchardef\nleqq="3\msb@14
\mathchardef\ngeqq="3\msb@15
\mathchardef\precneqq="3\msb@16
\mathchardef\succneqq="3\msb@17
\mathchardef\precnapprox="3\msb@18
\mathchardef\succnapprox="3\msb@19
\mathchardef\lnapprox="3\msb@1A
\mathchardef\gnapprox="3\msb@1B
\mathchardef\nsim="3\msb@1C
\mathchardef\napprox="3\msb@1D
\mathchardef\varsubsetneq="3\msb@20
\mathchardef\varsupsetneq="3\msb@21
\mathchardef\nsubseteqq="3\msb@22
\mathchardef\nsupseteqq="3\msb@23
\mathchardef\subsetneqq="3\msb@24
\mathchardef\supsetneqq="3\msb@25
\mathchardef\varsubsetneqq="3\msb@26
\mathchardef\varsupsetneqq="3\msb@27
\mathchardef\subsetneq="3\msb@28
\mathchardef\supsetneq="3\msb@29
\mathchardef\nsubseteq="3\msb@2A
\mathchardef\nsupseteq="3\msb@2B
\mathchardef\nparallel="3\msb@2C
\mathchardef\nmid="3\msb@2D
\mathchardef\nshortmid="3\msb@2E
\mathchardef\nshortparallel="3\msb@2F
\mathchardef\nvdash="3\msb@30
\mathchardef\nVdash="3\msb@31
\mathchardef\nvDash="3\msb@32
\mathchardef\nVDash="3\msb@33
\mathchardef\ntrianglerighteq="3\msb@34
\mathchardef\ntrianglelefteq="3\msb@35
\mathchardef\ntriangleleft="3\msb@36
\mathchardef\ntriangleright="3\msb@37
\mathchardef\nleftarrow="3\msb@38
\mathchardef\nrightarrow="3\msb@39
\mathchardef\nLeftarrow="3\msb@3A
\mathchardef\nRightarrow="3\msb@3B
\mathchardef\nLeftrightarrow="3\msb@3C
\mathchardef\nleftrightarrow="3\msb@3D
\mathchardef\divideontimes="2\msb@3E
\mathchardef\varnothing="0\msb@3F
\mathchardef\nexists="0\msb@40
\mathchardef\mho="0\msb@66
\mathchardef\thorn="0\msb@67
\mathchardef\beth="0\msb@69
\mathchardef\gimel="0\msb@6A
\mathchardef\daleth="0\msb@6B
\mathchardef\lessdot="3\msb@6C
\mathchardef\gtrdot="3\msb@6D
\mathchardef\ltimes="2\msb@6E
\mathchardef\rtimes="2\msb@6F
\mathchardef\shortmid="3\msb@70
\mathchardef\shortparallel="3\msb@71
\mathchardef\smallsetminus="2\msb@72
\mathchardef\thicksim="3\msb@73
\mathchardef\thickapprox="3\msb@74
\mathchardef\approxeq="3\msb@75
\mathchardef\succapprox="3\msb@76
\mathchardef\precapprox="3\msb@77
\mathchardef\curvearrowleft="3\msb@78
\mathchardef\curvearrowright="3\msb@79
\mathchardef\digamma="0\msb@7A
\mathchardef\varkappa="0\msb@7B
\mathchardef\hslash="0\msb@7D
\mathchardef\hbar="0\msb@7E
\mathchardef\backepsilon="3\msb@7F

\def\Bbb#1{\Bbb@#1}%
\def\Bbb@#1{{\msbm#1}}%
\def\frak#1{\frak@#1}%
\def\frak@#1{{\eul#1}}%

\makeatother

\makeatletter

\def\newremark#1{\@ifnextchar[{\@ormkwsq{#1}}{\@nrmkwsq{#1}}}

\def\@nrmkwsq#1#2{%
\@ifnextchar[{\@xnrmkwsq{#1}{#2}}{\@ynrmkwsq{#1}{#2}}}

\def\@xnrmkwsq#1#2[#3]{\expandafter\@ifdefinable\csname #1\endcsname
{\@definecounter{#1}\@addtoreset{#1}{#3}%
\expandafter\xdef\csname the#1\endcsname{\expandafter\noexpand
  \csname the#3\endcsname \@rmkcountersep \@rmkcounter{#1}}%
\global\@namedef{#1}{\@rmkwsq{#1}{#2}}
\global\@namedef{end#1}{\@endremarkwithsquare}}}

\def\@ynrmkwsq#1#2{\expandafter\@ifdefinable\csname #1\endcsname
{\@definecounter{#1}%
\expandafter\xdef\csname the#1\endcsname{\@rmkcounter{#1}}%
\global\@namedef{#1}{\@rmkwsq{#1}{#2}}
\global\@namedef{end#1}{\@endremarkwithsquare}}}

\def\@ormkwsq#1[#2]#3{\expandafter\@ifdefinable\csname #1\endcsname
  {\global\@namedef{the#1}{\@nameuse{the#2}}%
\global\@namedef{#1}{\@rmkwsq{#2}{#3}}%
\global\@namedef{end#1}{\@endremarkwithsquare}}}

\def\@rmkwsq#1#2{\refstepcounter
    {#1}\@ifnextchar[{\@yrmkwsq{#1}{#2}}{\@xrmkwsq{#1}{#2}}}

\def\@xrmkwsq#1#2{\@beginremark{#2}{\csname the#1\endcsname}\ignorespaces}
\def\@yrmkwsq#1#2[#3]{\@opargbeginremark{#2}{\csname
       the#1\endcsname}{#3}\ignorespaces}

\def\@rmkcounter#1{\noexpand\arabic{#1}}
\def\@rmkcountersep{.}
\def\@beginremark#1#2{\trivlist \item[\hskip \labelsep{\bf #1\ #2.}]}
\def\@opargbeginremark#1#2#3{\trivlist
      \item[\hskip \labelsep{\bf #1\ #2\ (#3)}]}
\def\@endremarkwithsquare{~\hspace{\fill}~$\square$\endtrivlist}
\makeatother

\makeatletter

   \def\@begintheorem#1#2{\sl \trivlist \item[\hskip \labelsep{\bf #2\ #1}]}
   \def\@opargbegintheorem#1#2#3{\sl \trivlist
            \item[\hskip \labelsep{\bf #1\ #2\ (#3)}]}

   \def\section{\@startsection {section}{1}{\z@}{-3.5ex plus -1ex minus
    -.2ex}{2.3ex plus .2ex}{\large\bf}} 
\makeatother

\makeatletter
\newenvironment{eqn}{\refstepcounter{subsection}
$$}{\leqno{\rm(\thesubsection)}$$\global\@ignoretrue}
\makeatother

\makeatletter
\newenvironment{subeqn}{\refstepcounter{subsubsection}
$$}{\leqno{\rm(\thesubsubsection)}$$\global\@ignoretrue}
\makeatother



\makeatletter
\def\@rmkcounter#1{\noexpand\arabic{#1}}
\def\@rmkcountersep{.}
\def\@beginremark#1#2{\trivlist \item[\hskip \labelsep{\bf #2\ #1.}]}
\def\@opargbeginremark#1#2#3{\trivlist
      \item[\hskip \labelsep{\bf #2\ #1\ (#3).}]}
\def\@endremarkwithsquare{~\hspace{\fill}~$\square$\endtrivlist}
\makeatother

\newenvironment{prf}[1]{\trivlist
\item[\hskip \labelsep{\it
#1.\hspace*{.3em}}]}{~\hspace{\fill}~$\square$\endtrivlist}

\newenvironment{proof}{\begin{prf}{\bf Proof}}{\end{prf}}

\let\tempcirc=\circ
\def\circ{\mathord{\raise0.25ex\hbox{$\scriptscriptstyle\tempcirc$}}}

\title{On the prime-to-$p$ part of the groups of connected components
of N\'eron models.}
\author{Bas Edixhoven}

\newtheorem{theorem}[subsection]{Theorem.}

\newtheorem{lemma}[subsection]{Lemma.}
\newtheorem{corollary}[subsection]{Corollary.}

\newtheorem{definition}[subsection]{Definition.}

\newremark{rmk}[subsection]{Remark}
\newremark{subrmk}[subsubsection]{Remark}
\newremark{notation}[subsection]{Notation}
\newremark{example}[subsection]{Example}

\renewcommand{\baselinestretch}{1.3}

\newcommand{\msy}{\bf}

\newcommand{\eps}{\varepsilon}
\newcommand{\Spec}{{\rm Spec}}
\newcommand{\Z}{{\msy Z}}
\newcommand{\Q}{{\msy Q}}
\newcommand{\C}{{\msy C}}
\newcommand{\R}{{\msy R}}
\newcommand{\ld}{\langle}
\newcommand{\rd}{\rangle}
\newcommand{\F}{{\msy F}}
\newcommand{\Gal}{{\rm Gal}}
\newcommand{\Pic}{{\rm Pic}}
\newcommand{\Gm}{{{\rm G}_m}}
\newcommand{\Ext}{{\rm Ext}}
\newcommand{\im}{{\rm im}}
\newcommand{\Ksep}{{K^{\rm s}}}
\newcommand{\Dsep}{{D^{\rm s}}}
\newcommand{\kbar}{{\overline{k}}}
\newcommand{\rT}{{\rm T}}
\newcommand{\rH}{{\rm H}}
\newcommand{\tors}{{\rm tors}}
\newcommand{\Itame}{{I_{\rm t}}}
\newcommand{\Vtilde}{{\widetilde{V}}}
\newcommand{\Wtilde}{{\widetilde{W}}}
\newcommand{\ttilde}{{\tilde{t}}}
\newcommand{\atilde}{{\tilde{a}}}
\newcommand{\Ktilde}{{\widetilde{K}}}
\newcommand{\Ktame}{{K^{\rm t}}}
\newcommand{\ttame}{{t_{\rm t}}}
\newcommand{\atame}{{a_{\rm t}}}
\newcommand{\utame}{{u_{\rm t}}}
\newcommand{\Dtame}{{D^{\rm t}}}
\newcommand{\rank}{{\rm rank}}
\newcommand{\Itilde}{{\widetilde{I}}}
\newcommand{\phitilde}{{\widetilde{\phi}}}
\newcommand{\phibar}{{\overline{\phi}}}
\newcommand{\ra}{{\rm a}}
\newcommand{\rt}{{\rm t}}

\begin{document}
\maketitle
\tableofcontents

\section{Introduction.}
The aim of this article was originally to improve certain results of Dino
Lorenzini concerning the groups of connected components of special fibres of
N\'eron models of abelian varieties. Let $D$ be a strictly henselian
discrete valuation ring, $K$ its field of fractions, $k$ its residue field and
$A_K$ an abelian variety over $K$ with N\'eron model $A$ over $D$. Let $p\geq0$
be the characteristic of $k$ and let $\Phi_{(p)}$ denote the prime-to-$p$ part
of the group of connected components $\Phi$ of $A_k$.
In \cite{Dino1} Lorenzini shows the existence of a functorial four step
filtration on $\Phi_{(p)}$ and he proves certain properties satisfied by
this filtration. In particular, he gives bounds for the successive quotients.
These bounds are of the following type. For a prime $l$ and a finite abelian
group $G$ of $l$-power order, say $G\cong\oplus_{i\geq1}\Z/l^{a_i}\Z$
with $a_1\geq a_2\geq\cdots$, he defines
$\delta_l'(G):=l^{a_1}-1+(l-1)\sum_{i\geq2}a_i$.
Then he gives bounds for the
$\delta_l'$ of certain successive quotients in terms of the dimensions of the
toric and abelian variety parts of the special fibres of N\'eron models of
$A_K$ over various extensions of~$D$. In \cite[Remark~2.16]{Dino1} he remarks
that the bounds might possibly be improved by replacing $\delta_l'$ by
an other invariant $\delta_l$ defined as follows: for $G$ as above one
has $\delta_l(G)=\sum_{i\geq1}(l^{a_i}-1)$. This improvement is exactly
what we do in this article. The results can be found in \S\ref{section3}.
Needless to say, we follow very much the approach of \cite{Dino1} in order
to prove these sharper bounds. In fact, only Lemma~2.13 of \cite{Dino1} has
to be changed, so the proof we give is rather short. We have taken this
opportunity to weaken slightly the hypotheses of Lorenzini's results
(he supposes $D$ to be complete and $k$ to be algebraically closed).

In \S\ref{section2} we recall Lorenzini's filtration. In \S\ref{section3} we
state and prove the bounds on the $\delta_l$ of the $l$-parts of certain
successive quotients and in \S\ref{section5} we show by some examples that
the bounds of \S\ref{section3} are sharp; \S\ref{section4} is used
to show some results on finite abelian groups that are needed in
the other sections.

After all this work it turned out that a complete classification of the
possible $\Phi_{(p)}$ for abelian varieties whose reduction has toric part,
abelian variety part and unipotent part of fixed dimensions was in reach.
The result, which is surprisingly simple to state, can be found in
Thm.~\ref{thm61}.

In this article we will frequently speak of the abelian variety part,
the toric part and the unipotent part of the fibre over $k$ of a N\'eron
model over $D$. Since we are only interested in the characteristic polynomials
of certain endomorphisms on the toric and abelian variety part, it suffices
to define these parts after base change to an algebraic closure of $k$, and
up to isogeny. Over the algebraic closure of $k$ we can apply Chevalley's
theorem; in \cite[Thm.~9.2.1, Thm.~9.2.2]{BLR} one finds what is needed, and
even more.

I would like to thank Xavier Xarles for indicating a mistake in an earlier
version of this article, and Rutger Noot for his help concerning the proofs
of Lemma's~\ref{lemma410} and~\ref{lemma411}.

\section{Lorenzini's filtration.} \label{section2}
Let $D$ be a discrete valuation ring, let $K$ be its field of fractions and
$k$ its residue field. Let $A_K$ be an abelian variety over $K$, $A$ its
N\'eron model over $D$ and $\Phi:=A_k/A_k^0$ the finite \'etale group scheme
over $k$ of connected components of the special fibre $A_k$. Let $p\geq0$ be
the characteristic of $k$ and let $\Phi_{(p)}$ be the prime-to-$p$ part of
$\Phi$; if $p=0$ we define $\Phi_{(p)}$ to be equal to $\Phi$.
In this section we will briefly recall the construction in \cite{Dino1} of a
descending filtration
\begin{eqn} \label{eqn21}
\Phi_{(p)} = \Phi_{(p)}^0 \supset \Phi_{(p)}^1 \supset \Phi_{(p)}^2 \supset
\Phi_{(p)}^3 \supset \Phi_{(p)}^4 = 0
\end{eqn}
which is functorial in $A_K$ and invariant under base change by automorphisms
of $D$.
Since $\Phi_{(p)}$ is the direct sum of its $l$-parts $\Phi_l$, with $l$
ranging through the primes different from $p$, it suffices to describe the
filtration on each $\Phi_l$. We replace $D$ by its strict henselization and
view the group scheme $\Phi$ over the separably closed field $k$ as just a
group.

Let $l\neq p$ be a prime number. Let $K\to\Ksep$ be a separable closure,
let $\Dsep$ be the integral closure of $D$ in $\Ksep$ and let $\kbar$ be
the residue field of $\Dsep$; note that $\kbar$ is an algebraic closure of
$k$ and that $k\to\kbar$ is purely inseparable. The first step in the
construction of the filtration is the description of $\Phi_l$ in
terms of the Tate module $U_l:=\rT_l(A(\Ksep))$ with its action by
$I:=\Gal(\Ksep/K)$ given in Prop.~11.2 of \cite{Grothendieck1}:
\begin{eqn} \label{eqn22}
\Phi_l = \left(U_l\otimes\Q/\Z\right)^I/
\left(U_l^I\otimes\Q/\Z\right)
\end{eqn}
The long exact cohomology sequence of the short exact sequence
$$
0\to U_l\to U_l\otimes\Q\to U_l\otimes\Q/\Z\to 0
$$
of continuous $I$-modules gives a canonical isomorphism
(see \cite[(11.3.8)]{Grothendieck1})
\begin{eqn} \label{eqn23}
\Phi_l = \tors(\rH^1(I,U_l))
\end{eqn}
where for $M$ any abelian group, $\tors(M)$ denotes the subgroup of
torsion elements. Let $\Itame$ be the quotient of $I$ corresponding to
the maximal tamely ramified extension of $D$, and let $P$ be the kernel
of $I\to\Itame$. Then $\Itame$ is canonically isomorphic to
$\prod_{q\neq p}\rT_q(\Gm(k))=\prod_{q\neq p}\Z_q(1)$ and $P$ is a pro-$p$
group. The Hochschild-Serre spectral sequence shows that
\begin{eqn}\label{eqn24}
\Phi_l = \tors(\rH^1(I,U_l)) = \tors(\rH^1(\Itame,U_l^P)) =
\tors((U_l^P)_\Itame)(-1) = \tors((U_l)_I)(-1)
\end{eqn}
with the lower indices $\Itame$ and $I$ denoting coinvariants and ``$(-1)$''
a Tate twist. Let $N_l$ be the submodule of $U_l$ which is generated
by the elements $\sigma(x)-x$ with $\sigma$ in $I$ and $x$ in $U_l$. Then
by definition we have $(U_l)_I=U_l/N_l$. As in \cite[\S2.5]{Grothendieck1}, we
define $V_l:=U_l^I$. Then $V_l$, which is called the fixed part of $U_l$, is
canonically isomorphic to $\rT_l(A_k(k))$.
Let $A_K'$ be the dual of $A_K$, i.e., $A_K'=\Pic^0_{A_K/K}$. We will denote
by $A'$ the N\'eron model over $D$ of $A_K'$, by $\Phi'$ its group of
connected components, etc.
Let $\ld{\cdot},{\cdot}\rd\colon U_l\times U_l'\to\Z_l(1)$ be the Weil
pairing.
For any $y$ in $V'_l$, $\sigma$ in $I$ and $x$ in $U_l$ we have
$\ld\sigma(x)-x,y\rd = \ld\sigma(x),y\rd-\ld x,y\rd =
\sigma(\ld x,\sigma^{-1}(y)\rd) - \ld x,y\rd = 0$. It follows that
$N_l$ is contained in the orthogonal ${V_l'}^\perp$ of $V_l'$ in $U_l$.
Since $U_l/{V_l'}^\perp$ is torsion free, we conclude that
\begin{eqn} \label{eqn25}
\Phi_l = \tors\left({V_l'}^\perp/N_l\right)(-1)
\end{eqn}
\begin{rmk} \label{rmk26}
In the proof of Thm.~\ref{thm33} we will see that ${V_l'}^\perp/N_l$ is in
fact a finite group, hence we have $\Phi_l = ({V_l'}^\perp/N_l)(-1)$.
\end{rmk}
Now it is clear that any filtration on ${V_l'}^\perp$ induces a filtration
on $\Phi_l$. As in \cite[\S2.5]{Grothendieck1}, we define $W_l\subset V_l$
to be the submodule corresponding to the maximal torus in $A_k$.
Let $\Wtilde_l\subset\Vtilde_l\subset V_l$ be the submodules called the
essentially toric part and the essentially fixed part in
\cite[\S4.1]{Grothendieck1}; if $G/k'$ is the connected component of the
special fibre of a semi-stable N\'eron model of $A_K$ over a suitable
sub-extension of $K\to\Ksep$ then $\Vtilde_l$ corresponds to $\rT_l(G(\kbar))$
and $\Wtilde_l$ to the Tate module of the maximal torus in $G$.
We denote by $t$, $a$ and $u$ the dimensions of the toric part, the
abelian variety part and the unipotent part of $A^0_\kbar$; we denote by
$\ttilde$ and $\atilde$ the analogous dimensions of any semi-stable reduction
of $A_K$. Note that $t+a+u=\ttilde+\atilde=\dim(A_K)$. An easy application
of the Igusa-Grothendieck orthogonality theorem (which states that
$W_l=V_l\cap{V_l'}^\perp$, see \cite[Thm.~2.4]{Grothendieck1}, or
\cite[Thm.~3.1]{Oort1}), gives us the following filtration of ${V_l'}^\perp$,
in which the successive quotients are torsion free and of the indicated rank:
\begin{eqn} \label{eqn27}
{V_l'}^\perp \;\;\stackrel{\ttilde-t}{\supset}\;\; \Vtilde_l\cap{V_l'}^\perp
\;\;\stackrel{2(\atilde-a)}{\supset}\;\; \Wtilde_l
\;\;\stackrel{\ttilde-t}{\supset}\;\; W_l
\;\;\stackrel{t}{\supset}\;\; 0
\end{eqn}
Lorenzini's filtration (\ref{eqn21}) on $\Phi_l$ is the filtration induced by
(\ref{eqn25}) and (\ref{eqn27}). Note that in fact any finite sub-extension
of $K\to\Ksep$ induces a filtration on $\Phi_l$ as above; see
\cite[Thm.~3.1]{Dino1} for results concerning those filtrations. The reason
we only consider the filtration coming from extensions over which $A_K$ has
semi-stable reduction is that only that filtration matters for the bounds
on $\Phi_{(p)}$ of the next section.

\section{Bounds on $\Phi_{(p)}$.} \label{section3}
We keep the notation of the previous section. Recall that $K$ is strictly
henselian. First we define some invariants of finite abelian groups and fix
some notation needed to state our results.
\begin{definition} \label{def31}
For $l$ a prime number and $a=(a_1,a_2,\ldots)$ a sequence of integers
$a_i\geq0$ with $a_i=0$ for $i$ big enough, let
$\delta_l(a):=\sum_i(l^{a_i}-1)$. For $l$ a prime number and
$G\cong\oplus_i\Z/l^{a_i}\Z$ a finite abelian group of $l$-power order let
$\delta_l(G):=\delta_l(a)$, where $a:=(a_1,a_2,\ldots)$.
For $G$ a finite abelian group let $\delta(G):=\sum_l\delta_l(G_l)$,
where $G=\oplus_l G_l$ is the decomposition of $G$ into groups of prime
power order.
\end{definition}

\begin{notation} \label{notation32}
Let $\Ktilde$ be the minimal sub-extension of $\Ksep$ over which $A_K$ has
semi-stable reduction; it corresponds to the kernel of $I$ acting on
$\Vtilde$, see \cite[\S4.1]{Grothendieck1}.
We define $\Ktame$ to be the maximal tame extension in $\Ktilde$, and
for all $l\neq p$ we let $K_l$ denote the maximal sub-extension of $\Ktilde$
whose degree over $K$ is a power of $l$. We denote by $\ttilde$, $\atilde$,
$\ttame$, $\atame$, $\utame$, $t_l$, $a_l$  and $u_l$ the dimensions of the
toric parts, the abelian variety parts and the unipotent parts
of the corresponding N\'eron models of $A_K$. For each prime $l\neq p$ we
let $I_{(l)}$ be the subgroup of $I$ such that $I/I_{(l)}$ is the quotient
$\Z_l(1)$ of $\Itame$.

Let $A^\rt$ be the N\'eron model of $A_K$ over the ring of integers $\Dtame$
of $\Ktame$. Then $\Gal(\Ktame/K)$ acts (from the right) on $A^\rt$,
compatibly with its right-action on $\Spec(\Dtame)$. This action induces an
action of $\Gal(\Ktame/K)$ on the special fibre $A^\rt_k$. Let $\sigma$ be
a generator of the cyclic group $\Gal(\Ktame/K)$. Let $l\neq p$ be a
prime number and $i\geq1$ an integer. Let $f_{l,i}$ denote the cyclotomic
polynomial whose roots are the roots of unity of order $l^i$.
We define $m_{\ra,l,i}$ and $m_{\rt,l,i}$ to be the multiplicities of
$f_{l,i}$ in the characteristic polynomials of $\sigma$ on the abelian variety
part and on the toric part, respectively, of $A^\rt_k$ (say one lets $\sigma$
act on $\rT_l(A^\rt_k(\kbar))\otimes\Q$).
Let $m_{l,i}:=m_{\ra,l,i}+m_{\rt,l,i}$. Finally, for $j\geq1$ we define
$p_{\ra,l,j}:=|\{i\geq1\;|\;m_{\ra,l,i}\geq j\}|$,
$p_{\rt,l,j}:=|\{i\geq1\;|\;m_{\rt,l,i}\geq j\}|$ and
$p_{l,j}:=|\{i\geq1\;|\;m_{l,i}\geq j\}|$. For an interpretation of
$p_{\ra,l}=(p_{\ra,l,1},p_{\ra,l,2},\ldots)$ in terms of
$m_{\ra,l}=(m_{\ra,l,1},m_{\ra,l,2},\ldots)$ etc. using partitions, see the
beginning of the proof of Lemma~\ref{lemma45}.
\end{notation}

\begin{theorem} \label{thm33}
Let $l\neq p$ be a prime number and consider the filtration (\ref{eqn21})
on $\Phi_l$. With the notations above, we have:
\begin{enumerate}
\item The group $\Phi_l^3$ can be generated by $t$ elements.
\item $\delta_l(\Phi_l^2/\Phi_l^3) \leq \delta_l(p_{\rt,l}) \leq t_l-t$.
\item $\delta_l(\Phi_l^1/\Phi_l^2) \leq \delta_l(p_{\ra,l}) \leq 2(a_l-a)$.
\item $\delta_l(\Phi_l/\Phi_l^1) \leq \delta_l(p_{\rt,l}) \leq t_l-t$.
\item $\delta_l(\Phi_l/\Phi_l^2) \leq \delta_l(p_l) \leq (t_l-t)+2(a_l-a)$.
\item $\delta_l(\Phi_l^1/\Phi_l^3) \leq \delta_l(p_l) \leq (t_l-t)+2(a_l-a)$.
\end{enumerate}
\end{theorem}

\begin{corollary} \label{cor34}
\begin{enumerate}
\item The group $\Phi_{(p)}^3$ can be generated by $t$ elements.
\item $\delta(\Phi_{(p)}^2/\Phi_{(p)}^3) \leq \ttame-t.$
\item $\delta(\Phi_{(p)}^1/\Phi_{(p)}^2) \leq 2(\atame-a).$
\item $\delta(\Phi_{(p)}/\Phi_{(p)}^1) \leq \ttame-t.$
\item $\delta(\Phi_{(p)}/\Phi_{(p)}^2) \leq (\ttame-t)+2(\atame-a).$
\item $\delta(\Phi_{(p)}^1/\Phi_{(p)}^3) \leq (\ttame-t)+2(\atame-a).$
\end{enumerate}
\end{corollary}

\begin{prf}{{\bf Proof} {\rm (of Thm.~\ref{thm33})}}
We begin with some generalities. We always have $M_I=(M_{I_{(l)}})_{\Z_l(1)}$.
The functors $M\mapsto M_{I_{(l)}}$ and $M\mapsto M^{I_{(l)}}$ are exact on
the category of finitely generated $\Z_l$-modules with continuous
$I_{(l)}$-action, and, for such modules, the canonical map
$M^{I_{(l)}}\to M_{I_{(l)}}$ is an isomorphism, hence $M_{I_{(l)}}$ is torsion
free if $M$ is torsion free. For $M$ a finitely generated $\Z_l$-module with
continuous $\Z_l(1)$-action we have $M_{\Z_l(1)}=M/(\sigma-1)M$ and
$M^{\Z_l(1)}=M[\sigma-1]$, where $\sigma$ is any topological generator
of $\Z_l(1)$.

Next we recall some general facts on the action of $I$ on $U_l$.
Let $\Itilde$ denote the subgroup $\Gal(\Ksep/\Ktilde)$ of $I$.
Then $\Itilde$ acts trivially on $\Vtilde_l=U_l^\Itilde$ and on
$U_l/\Vtilde_l$; the action of $\Itilde$ on $U_l$ factors through the biggest
pro-$l$ quotient $\Z_l(1)$ of $\Itilde$ and is given by an isogeny
$U_l/\Vtilde_l\to\Wtilde_l(-1)$ (see \cite[\S9.2, Thm.~10.4]{Grothendieck1}).
It follows that $N_l\cap\Wtilde_l$ is open in $\Wtilde_l$. The group $I$
acts on $\Vtilde_l$ via its finite quotient $\Gal(\Ktilde/K)=I/\Itilde$; this
action can be described in terms of an action of $I/\Itilde$ on the special
fibre of the N\'eron model of $A_\Ktilde$ (see \cite[\S4.2]{Grothendieck1}).
Dually, $I$ acts with finite image on $U_l/\Wtilde_l$.

As promised in Remark.~\ref{rmk26} we will show that $\Phi_l={V_l'}^\perp/N_l$.
It suffices to show that ${V_l'}^\perp$ and $N_l$ have the same rank.
We have $\rank(U_l/{V_l'}^\perp)=t+2a$. From the generalities at the beginning
of the proof it follows that
$U_l/N_l=(U_l)_I=((U_l)_{I_{(l)}})_{\Z_l(1)}=(U_l^{I_{(l)}})_{\Z_l(1)}$.
Let $\sigma$ be a topological generator of $\Z_l(1)$. The exact sequence
\begin{eqn} \label{eqn35}
0 \longrightarrow U_l^I \longrightarrow U_l^{I_{(l)}}
\;\;\stackrel{\sigma-1}{\longrightarrow} \;\; U_l^{I_{(l)}} \longrightarrow
\left(U_l^{I_{(l)}}\right)_{\Z_l(1)} \longrightarrow 0
\end{eqn}
shows that $\rank((U_l)_I)=\rank(U_l^I)=\rank(V_l)=t+2a$. In order to prove
Thm.~\ref{thm33} we may neglect the Tate twist in~(\ref{eqn25}).

By definition, we have $\Phi_l^3=W_l/N_l\cap W_l$. Since $W_l$
is a free $\Z_l$ module of rank~$t$, $\Phi_l^3$ can be generated by $t$
elements.

Let us now consider $\Phi_l^2/\Phi_l^3$.
Since $\Phi_l^2=\Wtilde_l/\Wtilde_l\cap N_l$, the group $\Phi_l^2/\Phi_l^3$ is
a quotient of $(\Wtilde_l/W_l)_I=((\Wtilde_l/W_l)_{I_{(l)}})_{\Z_l(1)}$.
Lemma~\ref{lemma44} implies that
$\delta_l(\Phi_l^2/\Phi_l^3)\leq
\delta_l(((\Wtilde_l/W_l)_{I_{(l)}})_{\Z_l(1)})$.
By the generalities above, $(\Wtilde_l/W_l)_{I_{(l)}}$ is isomorphic as
$\Z_l(1)$-module to $\Wtilde_l^{\Gal(\Ktilde/K_l)}/W_l$. Note that
$\Wtilde_l^{\Gal(\Ktilde/K_l)}$ is the Tate module of the toric part of the
special fibre of the N\'eron model of $A_K$ over the ring of integers of $K_l$,
and that $W_l=\Wtilde_l^{\Gal(\Ktilde/K)}$. It follows that for all $i\geq1$
the multiplicity of $f_{l,i}$ in the characteristic polynomial of a generator
$\sigma$ of $\Gal(\Ktame/K)$ on $(\Wtilde_l/W_l)_{I_{(l)}}$ is $m_{\rt,l,i}$
and that $1$ is not a root of this characteristic polynomial.
Applying Lemma~\ref{lemma45} and Cor.~\ref{cor46} gives the second part of the
theorem.

The proof of parts 3, 4, 5 and 6 of the theorem follows the same lines.
For example, $\Phi_l^1/\Phi_l^2$ is a quotient of
$((\Vtilde_l\cap{V_l'}^\perp)/\Wtilde_l)_I$. The group $I$ acts with finite
image on $(\Vtilde_l\cap{V_l'}^\perp)/\Wtilde_l$. The Grothendieck-Igusa
orthogonality theorem \cite[Thm.~2.4]{Grothendieck1} shows that
$(\Vtilde_l\cap{V'_l}^\perp/\Wtilde_l)_{I_{(l)}}$ has rank $2(a_l-a)$.
We have
\begin{eqn} \label{eqn36}
\left(\frac{\Vtilde_l\cap{V_l'}^\perp}{\Wtilde_l}\right)_I\otimes\Q =
\frac{(\Vtilde_l\cap{V_l'}^\perp)^I}{\Wtilde_l^I}\otimes\Q =
\frac{V_l\cap{V_l'}^\perp}{W_l}\otimes\Q =
\frac{W_l}{W_l}\otimes\Q = 0
\end{eqn}
which shows that the hypotheses of Lemma~\ref{lemma45} are satisfied.
Since $\Vtilde_l\cap{V_l'}^\perp/\Wtilde_l$ is isogenous to $\Vtilde_l/V_l$,
the multiplicities of the $f_{l,i}$ in the characteristic polynomial of a
generator $\sigma$ of $\Gal(\Ktame/K)$ on
$(\Vtilde_l\cap{V_l'}^\perp/\Wtilde_l)_{I_{(l})}$ are precisely the
$m_{\ra,l,i}$.

The proof of part 6 is entirely similar to the proofs of parts~2 and 3.
For parts~4 and 5 one notes that ${V'_l}^\perp/\Vtilde_l\cap {V'_l}^\perp$
is dual to $\Wtilde_l'/W_l'$, that ${V'_l}^\perp/\Wtilde_l$ is dual to
$\Vtilde_l'/V_l'$ and one uses that $A_K$ and $A_K'$ are isogenous.
\end{prf}

\begin{prf}{{\bf Proof} {\rm (of Cor.~\ref{cor34})}}
One just considers the factorization into irreducible factors of the
characteristic polynomial of a generator $\sigma$ of $\Gal(\Ktame/K)$ acting
on the semi-abelian variety part of~$A^\rt_k$.
\end{prf}

\section{Some abelian group theory.} \label{section4}
In this section we prove some results needed in the proof of Thm.~\ref{thm33}.
We fix a prime number $l$ and consider finite $\Z_l$-modules, i.e., finite
abelian groups of $l$-power order. Recall that there is a bijection between
the set of isomorphism classes of finite $\Z_l$-modules and the set of
partitions (i.e., sequences $m=(m_1,m_2,\ldots)$ of non-negative integers
such that $m_1\geq m_2\geq\cdots$ and $m_i=0$ for $i$ big enough):
a finite $\Z_l$-module $M$ corresponds to the partition $m=(m_1,m_2,\ldots)$
which satisfies $M\cong\oplus_{i\geq1}\Z/l^{m_i}\Z$.
To any partition $m$ we attach the number
$\delta_l(m):=\sum_{i\geq1}(l^{m_i}-1)$. Note that with these definitions,
we have $\delta_l(M)=\delta_l(m)$, with $\delta_l(M)$ as in Def.~\ref{def31}.

\begin{lemma} \label{lemma41}
Let $0\to B\to E\to A\to 0$ be an extension of finite $\Z_l$-modules. Let
$b=(b_1,b_2,\ldots)$, $e$ and $a$ denote their invariants.
Define $n_i:=a_i+b_i$ and $n=(n_1,n_2,\ldots)$. Let $m=(m_1,m_2,\ldots)$ be
the invariant of $A\oplus B$, i.e., $m$ is the sequence obtained by
reordering $(a_1,b_1,a_2,b_2,\ldots)$. Then we have $m\leq e\leq n$,
with ``$\leq$'' the lexicographical ordering.
\end{lemma}
\begin{proof}
Let us first prove that $e\leq n$. We use induction on $|E|$. We have
$e_1\leq a_1+b_1=n_1$ since $l^{a_1+b_1}$ kills $E$. If $e_1<n_1$ there is
nothing to prove, so we suppose that $e_1=n_1$. Choose any element $x$ in
$E$ of order $l^{e_1}$ and consider the subgroup it generates. We get a
diagram
\begin{subeqn} \label{eqn411}
\begin{array}{ccccccccc}
 & & 0 & & 0 & & 0 & & \\
 & & \downarrow & & \downarrow & & \downarrow & & \\
0 & \to & \Z/l^{b_1}\Z & \to & \Z/l^{e_1}\Z & \to & \Z/l^{a_1}\Z & \to & 0 \\
 & & \downarrow & & \downarrow & & \downarrow & & \\
0 & \to & B & \to & E & \to & A & \to & 0 \\
 & & \downarrow & & \downarrow & & \downarrow & & \\
0 & \to & B' & \to & E' & \to & A' & \to & 0 \\
 & & \downarrow & & \downarrow & & \downarrow & & \\
 & & 0 & & 0 & & 0 & &
\end{array}
\end{subeqn}
in which the rows and columns are exact. Now the columns are split, since
$l^{b_1}$ is the exponent of $B$, etc. Hence $b':=(b_2,b_3,\ldots)$,
$e':=(e_2,e_3,\ldots)$ and $a':=(a_2,a_3,\ldots)$ are the invariants of $B'$,
$E'$ and $A'$, respectively. The proof is finished by induction.

Let us now prove that $e\geq m$. By passing to Pontrjagin duals, if necessary,
we may assume that $b_1\geq a_1$. Then $m_1=b_1$. If $e_1>m_1$ there is
nothing to prove, hence we suppose that $e_1=m_1=b_1$. We choose any
element $x$ in $B$ of order $l^{b_1}$. Just as above we find a diagram
\begin{subeqn} \label{eqn412}
\begin{array}{ccccccccc}
 & & 0 & & 0 & &  & & \\
 & & \downarrow & & \downarrow & & & & \\
 & & \Z/l^{b_1}\Z & \stackrel{\rm id}{\to} & \Z/l^{b_1}\Z & &  & & \\
 & & \downarrow & & \downarrow & & & & \\
0 & \to & B & \to & E & \to & A & \to & 0 \\
 & & \downarrow & & \downarrow & & \downarrow & & \\
0 & \to & B' & \to & E' & \to & A & \to & 0 \\
 & & \downarrow & & \downarrow & & & & \\
 & & 0 & & 0 & & & &
\end{array}
\end{subeqn}
in which the columns are split. Induction finishes the proof.
\end{proof}

\begin{rmk} \label{rmk42}
It would be nice to have a complete description of the possible invariants
of extensions $E$ of finite $\Z_l$-modules $A$ by $B$ in terms of the
invariants of $A$ and $B$. In particular, are there more restrictions than
the following: those in Lemma~\ref{lemma41}, $e_i\geq a_i$ and $e_i\geq b_i$
for all~$i$, and the minimal number of generators for $E$ does not exceed the
sum of those numbers for $A$ and $B$? As Hendrik Lenstra pointed out to me,
the problem can be phrased in terms of Hall polynomials, see for example
\cite{MacDonald1}.
\end{rmk}

\begin{lemma} \label{lemma43}
Suppose that $a=(a_1,a_2,\ldots)$ and $b=(b_1,b_2,\ldots)$ are partitions
of $N$ (i.e., $\sum_{i\geq1}a_i=N=\sum_{i\geq1}b_i$) and that $a\geq b$
in the lexicographical ordering. Then $\delta_l(a)\geq\delta_l(b)$, with
equality if and only if $a=b$.
\end{lemma}
\begin{proof}
Consider the set $X$ of all partitions of $N$ with its lexicographical
ordering. From the inequality
$$
(l^{n+1}-1) + (l^{m-1}-1) > (l^n-1) + (l^m-1)
$$
satisfied for any integers $n\geq m$ it follows that $\delta_l\colon X\to\Z$
is strictly increasing.
\end{proof}

\begin{lemma} \label{lemma44}
\begin{enumerate}
\item For $M$ a finite $\Z_l$-module we have $\delta_l(M)\geq0$, with
equality if and only if $M=0$.
\item Let $0\to M'\to M\to M''\to 0$ be a short exact sequence of finite
$\Z_l$-modules. Then $\delta_l(M)\geq\delta_l(M')+\delta_l(M'')$, with
equality if and only if the sequence is split.
\item Let $0\to M'\to M\to M''\to 0$ be a short exact sequence of finite
$\Z_l$-modules. Suppose that $M$ is killed by $l^a$ and that $|M'|=l^{b}$.
Then $\delta_l(M)\leq \delta_l(M'')+b(l^a-l^{a-1})$.
\end{enumerate}
\end{lemma}
\begin{proof}
This follows directly from Lemmas~\ref{lemma41} and \ref{lemma43}.
\end{proof}

\begin{lemma} \label{lemma45}
Let $M$ be a finitely generated free $\Z_l$-module with an automorphism
$\sigma$ of finite order. Suppose that $M/(\sigma-1)M$ is finite, or,
equivalently, that the automorphism $\sigma\otimes1$ of the $\Q_l$-vector
space $M\otimes\Q$ does not have $1$ as eigenvalue. For $i\geq1$ let
$m_i$ be the multiplicity, in the characteristic polynomial of $\sigma$,
of the cyclotomic polynomial $f_i$ whose roots are the roots of unity of
order $l^i$. For each $j\geq1$, let $p_j:=|\{i\geq1\;|\; m_i\geq j\}|$.
Then $\delta_l(M/(\sigma-1)M)\leq \sum_{i\geq1}(l^{p_i}-1)$.
\end{lemma}
\begin{proof}
Let $q=(q_1,q_2,\ldots)$ be the partition obtained by reordering
$(m_1,m_2,\ldots)$. Then $p:=(p_1,p_2,\ldots)$ is what is usually called
the conjugate of $q$: when viewing partitions as Young diagrams, $p$ and
$q$ are obtained from each other by interchanging rows and columns.
In particular, we have $\sum_{i\geq1}p_i=\sum_{i\geq1}m_i$.

Let $n$ be the order of $\sigma$. Then $M$ becomes a module over the ring
$\Z_l[x]/(x^n-1)$. Let us write $n=l^rn'$ with $n'$ not divisible by $l$.
Then $\Z_l[x]/(x^n-1)$ is the product of the ring $\Z_l[x]/(x^{l^r}-1)$
by another ring $R$ and $x-1$ is invertible in $R$. This implies that
$M$ is the direct sum of two modules, one over $\Z_l[x]/(x^{l^r}-1)$ and
the other over $R$, and that the module over $R$ does not contribute
to $M/(\sigma-1)M$. Hence we have reduced the problem to the case where the
order of $\sigma$ is $l^r$.

Let $i_1<i_2<\cdots<i_{p_1}$ denote the integers $i\geq1$ such that
$m_i>0$. For $1\leq j\leq p_1$, let $F_j:=f_{i_j}$ be the corresponding
cyclotomic polynomials, and let $F:=F_1{\cdot}F_2\cdots F_{p_1}$.
Since $1$ is not an eigenvalue of $\sigma$ on $M\otimes\Q$, and $M$ is
torsion free as $\Z_l$-module, $M$ is a module over the ring
$A:=\Z_l[x]/(F)$. For any $A$-module $N$, we define $\overline{N}:=N/(x-1)N$.
Let us first note that for all $j$ we have $F_j(1)=l$. It follows that
$\overline{A}=\Z/l^{p_1}\Z$. For $N$ an $A$-module, $\overline{N}$ is
an $\overline{A}$-module, hence $l^{p_1}$ annihilates $\overline{N}$.

We claim that $|\overline{M}|=l^{\sum_{i\geq1}p_i}$. To prove this, note
that $|\overline{M}|=|\det(\sigma-1)|_l^{-1}$, with $|\cdot|_l$ the $l$-adic
absolute value on $\Q_l$, normalized by $|l|_l=1/l$. So in order to
compute $|\overline{M}|$ we may replace $M$ by any $\sigma$-stable lattice
$M'$ in $M\otimes\Q\cong\oplus_{i\geq1}(\Q_l[x]/(f_i))^{m_i}$. Taking
$M':=\oplus_{i\geq1}(\Z_l[x]/(f_i))^{m_i}$ and noting that $f_i(1)=1$ gives
the result.

Let $a=(a_1,a_2,\ldots)$ be the invariant of $\overline{M}$, i.e.,
$\overline{M}\cong\oplus_{i\geq1}\Z/l^{a_i}\Z$ and $a_1\geq a_2\geq\cdots$.
Note that $a$ and $p$ are partitions of the same number, hence in view of
Lemma~\ref{lemma43}, it suffices to show that $a\leq p$ in the lexicographical
ordering. Since $l^{p_1}$ annihilates $\overline{M}$, we have $a_1\leq p_1$.
If $a_1<p_1$ there is nothing to prove, so we assume that $a_1=p_1$. Let
$y$ be in $M$ such that its image $\overline{y}$ in $\overline{M}$ corresponds
to $(1,0,0,\ldots)$. Let $A'$ denote the submodule $Ay$ of $M$. Since
$M$ is free as a $\Z_l$-module, $A'$ is free as a $\Z_l$-module, and we
have $A'=\Z_l[x]/(G)$, with $G$ dividing $F$. Let $0\leq s\leq p_1$ be the
number of irreducible factors of $G$. Then we have $\overline{A'}=\Z/l^s\Z$.
We have a short exact sequence
\begin{subeqn} \label{eqn451}
0 \to A' \to M \to M' \to 0
\end{subeqn}
of $A$-modules, with $M'$ not necessarily free as $\Z_l$-module.
Multiplication by $x-1$ on this sequence induces an exact sequence
\begin{subeqn} \label{eqn452}
0 \to M'[x-1] \to \overline{A'} \to \overline{M} \to \overline{M'} \to 0
\end{subeqn}
The element $\overline{1}$ of $\overline{A'}$, which is annihilated by $l^s$,
is mapped to $\overline{y}$ which has annihilator $l^{p_1}$. It follows that
$s=p_1$, that $G=F$ and that $M'[x-1]=0$. Let us now consider the finite
$A$-module $\tors(M')$. Multiplication by $x-1$ acts injectively, hence
bijectively. Since $x-1$ is in the maximal ideal of $A$, it follows that
$\tors(M')=0$, hence that $M'$ is free as $\Z_l$-module. The proof is now
finished by induction on $\rank(M)$, since
$\overline{M'}\cong\oplus_{i\geq2}\Z/l^{a_i}\Z$ and the partition $p'$
obtained from $M'$ is $(p_2,p_3,\ldots)$.
\end{proof}

\begin{rmk}\label{rmk46}
Lemma~\ref{lemma45} can be seen as a bound on the cohomology group
$\rH^1(\Z/n\Z,M)$, where $1$ in $\Z/n\Z$ acts on $M$ via $\sigma$. It is
an interesting question, raised by Xavier Xarles, to obtain similar
bounds for non-cyclic groups.
\end{rmk}

\begin{corollary} \label{cor46}
Let $M$ be a finitely generated free $\Z_l$-module with an automorphism
$\sigma$ of finite order. Suppose that $M/(\sigma-1)M$ is finite.
Then $\delta_l(M/(\sigma-1)M)\leq\rank(M)$.
\end{corollary}
\begin{proof}
We use the notation of the beginning of the proof of Lemma~\ref{lemma45}.
Then one has:
\begin{subeqn} \label{eqn461}
\rank(M)\geq \sum_{i\geq1}m_i\phi(l^i) \geq \sum_{i\geq1}q_i\phi(l^i)
= \sum_{i\geq1}\sum_{j=1}^{p_i}\phi(l^j) = \sum_{i\geq1}(l^{p_i}-1)
\end{subeqn}
The proof is finished by applying Lemma~\ref{lemma45}.
\end{proof}

\begin{lemma} \label{lemma47}
Let $M$ be a finitely generated free $\Z_l$-module with an automorphism
$\sigma$ of finite order. Suppose that $M/(\sigma-1)M$ is finite and that
$\delta_l(M/(\sigma-1)M)=\rank(M)$. Then $M$ is a direct sum of
$\Z_l$-modules of the type $\Z_l[x]/(f_1{\cdot}f_2\cdots f_r)$ with
$\sigma$ acting as multiplication by $x$ and where $f_i$ denotes the
cyclotomic polynomial whose roots are the roots of unity of order $l^i$.
\end{lemma}
\begin{proof}
The proof is by induction on $\rank(M)$ and consists of an inspection of the
proofs of Lemma~\ref{lemma45} and Cor.~\ref{cor46}. First of all we must have
that $n'=1$. Secondly, we note that $m_i=q_i$ for all $i\geq1$ since the
inequalities in (\ref{eqn461}) are equalities (here we use that
$\sum_{i\geq1}m_i=\sum_{i\geq1}q_i$ and that $\phi(l^i)<\phi(l^j)$ if
$1\leq i<j$). So $m$ is the conjugate partition of $p$, hence
$A=\Z_l[x]/(F)$ with $F=f_1{\cdot}f_2\cdots f_{p_1}$. The formula for the
number of elements of $|\overline{M}|$ in the proof of Lemma~\ref{lemma45}
shows that $a$ and $p$ are partitions of the same number. By the hypotheses
of the lemma we are proving, we have $\delta_l(a)=\delta_l(p)$.
Lemma~\ref{lemma43} implies that $a=p$. The end of the proof of
Lemma~\ref{lemma45} shows that $A'=A$ and that $M'$ is free as $\Z_l$-module.
By induction on $\rank(M)$, we know that $M'$ is of the indicated type.
It remains to show that the short exact sequence
(\ref{eqn451}) splits. To do that, it is sufficient to show that
$\Ext^1_A(A_i,A)=0$, where $A_i=\Z_l[x]/(f_1\cdots f_i)$ with $i\leq p_1$.
This $\Ext^1$ is easily computed using the projective
resolution
\begin{subeqn}\label{eqn471}
\cdots \longrightarrow A \;\;\stackrel{f}{\longrightarrow}\;\;A
\;\;\stackrel{g}{\longrightarrow}\;\;A
\;\;\stackrel{f}{\longrightarrow}\;\;A \longrightarrow A_i \longrightarrow 0
\end{subeqn}
with $f=f_1\cdots f_i$ and $g=f_{i+1}\cdots f_{p_1}$.
\end{proof}
The following lemmas will be used in \S\ref{section5} and \S\ref{section6}.
\begin{lemma} \label{lemma48}
Let $M$ be a finite $\Z_l$-module and let
$$
M=M^0\supset M^1\supset M^2\supset\cdots\supset M^r=0
$$
be a strictly descending filtration. Suppose that for all $i$ with
$0\leq i\leq r-2$ the group $M^i/M^{i+2}$ is cyclic. Then $M$ is cyclic.
\end{lemma}
\begin{proof}
For $r\leq2$ there is nothing to prove. If we know the result for $r=3$, the
general case follows by induction since then $M^0/M^3$ is cyclic and
the filtration $M^0\supset M^2\supset M^3\supset\cdots\supset M^r$ has length
$r-1$. So assume now that $r=3$. Let $x$ be an element of $M^0$ such that its
image in $M^0/M^2$ is a generator. Then a certain multiple $ax$ of $x$ gives a
generator of $M^1/M^2$. Since $M^1$ is cyclic, and $M^1/M^2$ a non-trivial
quotient, $ax$ is a generator of $M^1$. The subgroup of $M^0$ generated by $x$
contains $M^1$ and its quotient by $M^1$ is $M^0/M^1$. We conclude that $x$
generates~$M^0$.
\end{proof}

\begin{rmk}\label{remark49}
The proof of Lemma~\ref{lemma48} generalizes immediately to a proof of the
following assertion: let $A$ be a local ring and $M$ an $A$-module with
a finite strictly descending filtration $M^i$ such that the $M^i/M^{i+2}$ are
cyclic, then $M$ is cyclic.
\end{rmk}

\begin{lemma}\label{lemma410}
Let $0\to B\to E\to A\to 0$ be a short exact sequence of finite $\Z_l$-modules
with invariants $b$, $e$ and $a$. Let $t\geq0$ be an integer and suppose that
$B$ is generated by $t$ elements. Then for all $i\geq1$ we have
$a_i\geq e_{i+t}$.
\end{lemma}
\begin{proof}
For a partition $p$, let $p'$ denote its conjugate. Then for all $i\geq1$
we have $l^{a_i'}=|A[l^i]A[l^{i-1}]|$. Let $d$ be the endomorphism of the
set of partitions defined by: $d(p)_i=p_{i+1}$ for all $i\geq1$. Let $d'$
be the conjugate of $d$: $d'(p)=d(p')'$. Then $d'(p)_i=\max(0,p_i-1)$.
When viewing a partition $p$ as a Young diagram in which the $p_i$ are the
lengths of the columns, $d$ and $d'$ remove the longest column and row,
respectively. Note that for a finite $\Z_l$-module $M$ with invariant $m$,
the submodule $lM$ has invariant $d'(m)$. Note that $d$ and $d'$ commute.
In the rest of this proof we will consider the partial
ordering on the set of partitions in which $p\leq q$ if and only if for all
$i\geq1$: $p_i\leq q_i$. Note that $p\leq q$ is equivalent to $p'\leq q'$.
Below we will use that $p\geq q$ if and only if: $p_1'\geq q_1'$ and
$d'(p)\geq d'(q)$. We will also use that if $N$ and $M$ are finite
$\Z_l$-modules with invariants $n$ and $m$ such that $N$ is a subquotient of
$M$, then $n\leq m$.

The proof of the lemma is by induction on $|E|$. What we have to prove is that
$a\geq d^t(e)$. The exact sequence $0\to B[l]\to E[l]\to A[l]$ shows that
$a_1'\geq e_1'-t$. Note that $d^t(e)'_1=\max(0,e_1'-t)$, hence we have
$a_1'\geq d^t(e)'_1$. The exact sequence
$0\mapsto B\cap lE\to lE\to lA$ shows (induction hypothesis) that
$d'(a)\geq d^t(d'(e))=d'(d^t(e))$. The two inequalities just proved imply
that $a\geq d^t(e)$.
\end{proof}

\begin{lemma}\label{lemma411}
Let $l$ be a prime. Let $0\to B\to E\to A$ be a short exact sequence of
finite $\Z_l$-modules with invariants $b$, $e$ and $a$, respectively.
Then
\begin{subeqn}\label{eqn4111}
\delta_l(b) + \delta_l(a) \geq
\sum_{i\geq1}\left(
\frac{l^{\lfloor e_i/2\rfloor}+l^{\lceil e_i/2\rceil}}{2}-1\right)
\end{subeqn}
where for any real number $x$, $\lfloor x\rfloor$ and $\lceil x\rceil$ denote
the largest (resp. smallest) integer $\leq x$ (resp. $\geq x$).
\end{lemma}
\begin{proof}
We have $\sum_i(a_i+b_i)=\sum_i e_i$. Lemma~\ref{lemma41} asserts that
$a+b\geq e$ in the lexicographical ordering. Consider the set $S$ of all
pairs $(r,s)$ of partitions, such that $r+s\geq e$ and
$\sum_i(r_i+s_i)=\sum_i e_i$. Let $f\colon S\to \Z$ be the map which sends
$(r,s)$ to $\delta_l(r)+\delta_l(s)$. We will show that $f$ achieves its
minimum at all $(r,s)$ in $S$ with the property that, for all $i\geq1$,
one has $\{r_i,s_i\}=\{\lfloor e_i/2\rfloor,\lceil e_i/2\rceil\}$.

Suppose now that $(r,s)$ is an element of $S$ where $f$ has a minimum.
We have to show that $|r_i-s_i|\leq1$ for all $i\geq1$. Suppose that this is
not the case. Let $j\geq1$ be minimal for the property that $|r_j-s_j|>1$
and $|r_i-s_i|\leq1$ for all $i<j$. We may and do suppose that $r_j-s_j>1$.
Note that if $j>1$ we have $s_{j-1}>s_j$. We define $r'$ and $s'$ as
follows: $(r'_i,s'_i)=(r_i-1,s_i+1)$ if $i\geq j$ and $r_i=r_j$;
in all other cases $(r'_i,s'_i)=(r_i,s_i)$. Note that $r'$ and $s'$ are
partitions, that $\sum_i(r'_i+s'_i)$ is equal to $\sum_i e_i$ and that
$f(r',s')$ is strictly smaller than $f(r,s)$.
\end{proof}

\section{Examples.} \label{section5}
The aim of this section is to give examples that show that the bounds in
Thm.~\ref{thm33} and Cor.~\ref{cor34} are sharp, in a sense that will become
clear in the examples. The examples we construct here will play an important
role in \S\ref{section6}.
We give our examples over the field $K:=\C((q))$ of formal Laurent series
over the complex numbers with its usual valuation, but it is easy to get
similar examples in mixed characteristic, or equal characteristic $p>0$.

The building stones of our examples are the following. For each integer
$n\geq1$ we let $E_n$ be the so-called Tate elliptic curve ``$\Gm/q^{n\Z}$''
over $K$ as described in \cite[\S VII]{DeligneRapoport} or in
\cite[\S6]{Mumford} ($E_n$ is obtained from the analytic family of elliptic
curves over the punctured unit disc with coordinate $q$ whose fibres are the
$\C^*/q^{n\Z}$, by base change from the field of finite tailed convergent
Laurent series to $K$). It is well known that the special fibre of the N\'eron
model of $E_n$ over $D:=\C[[q]]$ is an extension of $\Z/n\Z$ by the
multiplicative group.

For each prime $l$ and integer $r\geq0$ we define the ring
$\Lambda_{l,r}:=\Z[x]/(f_{l,1}\cdots f_{l,r})$, where as before $f_{l,i}$ is
the polynomial whose roots are the roots of unity of order $l^i$.
When $l>2$, we let $A_{l,r}$ be an abelian variety over $\C$ obtained as
follows: we choose an isomorphism of $\R$-algebras between
$\Lambda_{l,r}\otimes\R$ and a product of a number of copies of $\C$ and
define $A_{l,r}:=(\Lambda_{l,r}\otimes\R)/\Lambda_{l,r}$ (it is well known
that the trace form on $\Lambda_{l,r}$ implies the existence of a
polarization). The first three examples will be isogenous to twists of
products of copies of $E_{n}$ and of $A_{l,r,K}$. Of course
Lemma~\ref{lemma47} tells us how to cook up the required examples.

\begin{example}\label{example51}
Let $d\geq0$ and let $G$ be any finite abelian group that can be generated by
$d$ elements. Then $G\cong\oplus_{i=1}^d\Z/n_i\Z$, say.
For $A_K:=\prod_{i=1}^d E_{n_i}$ one has $\Phi=\Phi^3=\oplus_{i=1}^t\Z/n_i\Z$
and we have $d=t=\dim(A_K)$.
\end{example}

\begin{example}\label{example52}
Now consider parts~2 and 4 of Thm.~\ref{thm33}. Let $l$ be a prime.
For $i\geq1$ we let $B_{l,i}$ be the abelian variety $E_1\otimes\Lambda_{l,i}$
over $\C$, i.e., $B_{l,i}$ is a direct sum of copies of $E_1$, indexed by
some $\Z$-basis of $\Lambda_{l,i}$, and $\Lambda_{l,i}$ acts on $B_{l,i}$
according to its action on itself. In particular, multiplication by $x$ in
$\Lambda_{l,i}$ induces an automorphism $\sigma$ of $B_{l,i}$. Note that
$\sigma$ has order $l^i$. Let $C_{l,i}$ be the twist of $B_{l,i,K}$ over
$K(q^{1/l^i})$ by $\sigma$, i.e., $C_{l,i}$ is the quotient of the $K$-scheme
$B_{l,i,K}\times_{\Spec(K)}\Spec(K(q^{1/l^i}))$ by the group
$\Gal(K(q^{1/l^i})/K)=\Z/l^i\Z$ (here we choose a root of unity of order $l^i$)
which acts by $a\mapsto\sigma^a$ on the first
factor and via its natural action on the second factor.

We will now compute the group of connected
components $\Psi$ of the N\'eron model of $C_{l,i}$ over $D$, using
(\ref{eqn24}).
First of all we have $\rT_l(E_1(\Ksep))=\Z_l(1)\oplus\Z_l$, with $I$
acting via its quotient $\Z_l(1)$ in the following way: an element of $I$ with
image $a$ in $\Z_l(1)$ acts as multiplication by the matrix
$({1\atop0}{a\atop1})$. By construction,
$\rT_l(C_{l,i}(\Ksep))=\rT_l(E_1(\Ksep))\otimes\Lambda_{l,i}$
and an element in $I$ with image $a$ in $\Z_l(1)$ acts as
$({1\atop0}{a\atop1})\otimes x^a$. Since $C_{l,i}$ has $\atilde=t=0$, we have
$\Psi_l^1=\Psi_l^2$ and $\Psi_l^3=0$. The filtration
$\Z_l(1)\subset\rT_l(E_1(\Ksep)$ induces the filtration
$\Wtilde_l\subset\rT_l(C_{l,i}(\Ksep))$. It follows that $\Psi_l$ is the
cokernel of $({x-1\atop0}{x\atop x-1})$ and that $\Psi_l/\Psi_l^1$ and
$\Psi_l^2/\Psi_l^3$ are both isomorphic to $\Lambda_{l,i}/(x-1)=\Z/l^i\Z$.
An analogous computation shows that $\Psi=\Psi_l$. One can show that
$\Psi$ is isomorphic to $\Z/l^i\Z\oplus\Z/l^i\Z$ if $l>2$ and to
$\Z/2^{i+1}\Z\oplus\Z/2^{i-1}\Z$ if $l=2$.

Let $G$ be a finite abelian group of $l$-power order, say with invariant
$a=(a_1,a_2,\ldots)$. Then for $A_K:=\prod_{i\geq1}C_{l,a_i}$ we have
$\Phi_l/\Phi_l^1\cong\Phi_l^2/\Phi_l^3\cong G$ and
$\delta_l(G)=t_l=\dim(A_K)$. We remark that abelian varieties over $K$ that
are isogeneous to $A_K$ provide examples with $\Phi_l/\Phi_l^1$ not isomorphic
to $\Phi_l^2/\Phi_l^3$.
\end{example}

\begin{example}\label{example53}
For $l>2$ prime and $i\geq0$ we let $D_{l,i}$ be the abelian
variety over $K$ obtained by twisting $A_{l,i,K}$ over $K(q^{1/l^i})$ by the
automorphism $\sigma$ of $A_{l,i}$ which is induced from the multiplication
by $x$ in $\Lambda_{l,i}$. Then we have
$\rT_l(D_{l,i}(\Ksep))=\rT_l(A_{l,i}(\C))=\Lambda_{l,i}\otimes\Z_l$, and an
element in $I$ with image $a$ in $\Z_l(1)$ acts as $x^a$. Let $\Psi_l$
denote the group of connected components of attached to $D_{l,i}$. In this
case we have $\ttilde=a=0$, hence $\Psi_l=\Psi_l^1$ and $\Psi_l^2=0$.
By (\ref{eqn24}) we have $\Psi_l=\Lambda_{l,i}/(x-1)=\Z/l^i\Z$.

Suppose now that $l\neq2$. Let $G$ be a finite abelian group of $l$-power
order, say with invariant $a=(a_1,a_2,\ldots)$. Then for
$A_K:=\prod_{i\geq1}D_{l,a_i}$ we have
$\Phi_l=\Phi_l^1$, $\Phi_l^2=0$, $\Phi_l^1/\Phi_l^2\cong G$ and
$\delta_l(G)=2a_l=2\dim(A_K)$.
The case $l=2$ is a little bit different because $m_{\ra,2,1}$ is always even.
\end{example}

\begin{example} \label{example54}
Let $l$ be prime and let $r>0$ and $s>0$ be positive integers.
We will construct an abelian variety $A_K$ with $t=a=0$,
$\ttilde=l^r-1$, $\atilde=(l^{r+s}-l^r)/2$ and $\Phi=\Phi_l\cong\Z/l^{2r+s}\Z$.
It follows from Thm.~\ref{thm33} that in such an example $\Phi_l/\Phi_l^1$ and
$\Phi_l^2$ are cyclic of order $l^r$, that $\Phi_l^1/\Phi_l^2$ is
cyclic of order $l^s$ and that $\Phi_l/\Phi_l^2$ and $\Phi_l^1$ are cyclic of
order $l^{r+s}$. Hence this example shows that, as far as the exponent is
concerned, the two-fold extension $\Phi_l/\Phi_l^3$ can be arbitrary.

As in the previous examples, $f_{l,i}$ will denote the polynomial whose
roots are the roots of unity of order $l^i$, and $\Lambda_{l,r}$ is the
ring $\Z[x]/(f_{l,1}\cdots f_{l,r})$.
Let $\Lambda_{l,r,s}:=\Z[x]/(f_{l,r+1}\cdots f_{l,r+s})$. Let $D_{l,r,s}$ be
an abelian variety over $K$ obtained by replacing $\Lambda_{l,i}$ by
$\Lambda_{l,r,s}$ and $q^{1/l^i}$ by $q^{1/l^{r+s}}$ in the construction of
$D_{l,i}$ in Example~\ref{example53}.
Let $C_{l,r}$ be as in Example~\ref{example52}. Our example $A_K$ will be
isogeneous to $C_{l,r}\times D_{l,r,s}$.
Let $V:=\rT_l((C_{l,r}\times D_{l,r,s})(\Ksep))\otimes\Q$. Then $V$ is a
$\Q_l$-vector space with an action of $I=\Gal(\Ksep/K)$. We have an
isomorphism of $\Q_l$-vector spaces with $I$-action
\begin{subeqn}\label{eqn541}
\Lambda_{l,r}\otimes\Q_l \;\oplus\; \Lambda_{l,r,s}\otimes\Q_l \;\oplus\;
\Lambda_{l,r}\otimes\Q_l \;\;\tilde{\longrightarrow}\;\; V
\end{subeqn}
such that an element of $I$ with image $a$ in $\Z_l(1)$ acts via
\begin{subeqn}\label{eqn542}
\left(\begin{array}{ccc}x^a&0&ax^a\\0&x^a&0\\0&0&x^a\end{array}\right)
\end{subeqn}
Let
\begin{subeqn}\label{eqn543}
V=V^0\supset V^1\supset V^2\supset V^3=0
\end{subeqn}
be the filtration (\ref{eqn27}) on $V$. Then $V^2$ is simply the first
term in (\ref{eqn541}) and $V^1$ is the sum of the first two terms.
For any $\Z_l$-lattice $M$ in $V$ let $M^i:=M\cap V^i$.
To get our example $A_K$, it suffices to find an $I$-invariant $\Z_l$-lattice
$M$ in $V$ such that $M^1$ and $M/M^2$ are isomorphic, as $\Z_l[I]$-modules,
to $\Lambda_{l,r+s}\otimes\Z_l$, where an element of $I$ with image $a$ in
$\Z_l(1)$ acts on $\Lambda_{l,r+s}\otimes\Z_l$ as $x^a$. Namely, since $M$ is
$I$-invariant, $M$ is the $l$-adic Tate module of an abelian variety $A_K$
which is isogeneous to $C_{l,r}\times D_{l,r,s}$; for $A_K$ one has
$\Phi_l/\Phi_l^2$ and $\Phi_l^1$ cyclic of order $l^{r+s}$, hence $\Phi_l$
cyclic of order $l^{2r+s}$ by Lemma~\ref{lemma48}.

Let us now try to find such a $M$. Note that we have canonical projections
$\Lambda_{l,r+s}\to\Lambda_{l,r}$ and $\Lambda_{l,r+s}\to\Lambda_{l,r,s}$
which induce an embedding $\Lambda_{l,r+s}\otimes\Z_l\subset V^1$.
We will first show that we only have to
look among the sublattices $M$ with $M^1=\Lambda_{l,r+s}\otimes\Z_l$.
Namely, if $M$ is an $I$-invariant $\Z_l$-lattice in $V$ of the type we are
looking for, then for a suitable element of the form $v=(a,b,a)$ of $V^*$
(here we consider $V$ as a $\Q_l$-algebra and $V^*$ denotes the group of units
of $V$) $vM$ is isomorphic to $M$ as $\Z_l[I]$-module and has
$(vM)^1=\Lambda_{l,r+s}\otimes\Z_l$. Since the $\Z_l[I]$-module structure
determines the filtration, we also have $(vM)/(vM)^2\cong M/M^2$.
\par From now on we only consider $M$ with $M^1=\Lambda_{l,r+s}\otimes\Z_l$.
Such $M$ are determined by their image in $V/M^1$. So we look for an
$I$-invariant torsion free $\Z_l$-submodule $N$ of
$V/M^1=V^1/M^1\oplus\Lambda_{l,r}\otimes\Q_l$ whose image in
$\Lambda_{l,r}\otimes\Q_l$ is a lattice and for whose associated $M$ we have
$M/M^2\cong\Lambda_{l,r+s}\otimes\Z_l$. It follows that such a $N$ is
isomorphic, via the canonical projection, to its image in
$\Lambda_{l,r}\otimes\Q_l$. Lemma~\ref{lemma47} implies that this image is
isomorphic to $\Lambda_{l,r}\otimes\Z_l$, hence of the form
$z{\cdot}\Lambda_{l,r}\otimes\Z_l$ for some $z$ in
$(\Lambda_{l,r}\otimes\Q_l)^*$. We conclude that $N$ is of the form
$\im(\alpha)$, where
\begin{subeqn}\label{eqn544}
\alpha\colon \Lambda_{l,r}\otimes\Z_l \longrightarrow
V^1/M^1\;\oplus\;\Lambda_{l,r}\otimes\Q_l, \quad
a\mapsto(\phi(a),za)
\end{subeqn}
with $\phi\colon\Lambda_{l,r}\otimes\Z_l\to V^1/M^1$ a morphism of
$\Z_l$-modules, and $z\in(\Lambda_{l,r}\otimes\Q_l)^*$. For a given pair
$(\phi,z)$, let $N_{\phi,z}$ denote the image of the corresponding $\alpha$.

Let us first study what it means for $(\phi,z)$ that $N_{\phi,z}$ is
$I$-invariant. Using that $N_{\phi,z}$ is $I$-invariant if and only if
it is invariant under the matrix in (\ref{eqn542}) with $a$ replaced by $1$,
one easily sees that $N_{\phi,z}$ is $I$-invariant if and only if
\begin{subeqn}\label{eqn545}
\forall a\in\Lambda_{l,r}\otimes\Z_l\colon \quad
\phi(xa)=x\phi(a)+\overline{(xza,0)}
\end{subeqn}
To find out which $(\phi,z)$ satisfy (\ref{eqn545}), we write out everything
in terms of the $\Z_l$-basis $(1,x,\ldots,x^{l^r-2})$ of
$\Lambda_{l,r}\otimes\Z_l$. Let $y=(y_1,y_2)$ be in $V^1$ such that
$\phi(1)=\overline{y}$. One then checks that
\begin{subeqn}\label{eqn546}
\phi(x^i) = x^i\overline{y} + ix^i\,\overline{(z,0)}, \quad
\mbox{for $0\leq i\leq l^r-2.$}
\end{subeqn}
Applying (\ref{eqn545}) with $a=x^{l^r-2}$, and using that
$\sum_{i=0}^{l^r-1}x^i=0$ in $\Lambda_{l,r}$, gives
\begin{subeqn}\label{eqn547}
(xg_r'(x)z,g_r(x)y_2) \in \Lambda_{l,r+s}\otimes\Z_l
\end{subeqn}
where $g_r=f_{l,1}\cdots f_{l,r}$ and $g_r'$ is the derivative of $g_r$. The
conclusion is that $N_{\phi,z}$ is $I$-invariant if and only if $\phi$ is
given by (\ref{eqn546}) and $(y,z)$ satisfies (\ref{eqn547}). For a given
such pair $(y,z)$, let $M_{y,z}$ denote the lattice $M$ in $V$ corresponding
to $N_{\phi,z}$.

It remains now to be seen that there exist $(y,z)$ satisfying (\ref{eqn547})
such that $M_{y,z}/M_{y,z}^2$ is isomorphic to $\Lambda_{l,r+s}\otimes\Z_l$,
or, equivalently, such that $(M_{y,z}/M_{y,z}^2)_I$ is cyclic.
In order to have a useful description of $M_{y,z}$, we lift
$\phi$ to $V^1$ as follows: let
$\phitilde\colon\Lambda_{l,r}\otimes\Z_l\to V^1$ be the morphism of
$\Z_l$-modules such that
\begin{subeqn}\label{eqn548}
\phitilde\colon x^i\mapsto x^iy + ix^i(z,0),
\quad\mbox{for $0\leq i\leq l^r-2$}
\end{subeqn}
Then we have an isomorphism of $\Z_l$-modules:
\begin{subeqn}\label{eqn549}
\beta\colon \Lambda_{l,r+s}\otimes\Z_l\;\oplus\;\Lambda_{l,r}\otimes\Z_l
\;\;\tilde{\longrightarrow}\;\; M\subset V,\quad
(a,b)\mapsto (a,0)+(\phitilde(b),zb)
\end{subeqn}
Let $\tau$ be an element of $I$ with image $1$ in $\Z_l(1)$. Then $\tau$ acts
on $V$ by the matrix in (\ref{eqn542}) with $a=1$. One computes that in order
to make $\beta$ invariant under $I$, one must let $\tau$ act on the source of
$\beta$ in (\ref{eqn549}) by
\begin{subeqn}\label{eqn5410}
\tau\colon (a,b)\mapsto \left(xa+(xzb,0)+x\phitilde(b)-\phitilde(xb),xb\right)
\end{subeqn}
Using this formula, we can study $M_{y,z}/M_{y,z}^2$. Recall that
$\Lambda_{l,r+s}\otimes\Z_l$ is the image in $V^1$ of the sum of the two
canonical projections from $\Lambda_{l,r+s}\otimes\Z_l$ to
$\Lambda_{l,r}\otimes\Z_l$ and $\Lambda_{l,r,s}\otimes\Z_l$. It follows that
\begin{subeqn}\label{eqn5411}
M_{y,z}/M_{y,z}^2 \cong N:=
\Lambda_{l,r,s}\otimes\Z_l \;\oplus\; \Lambda_{l,r}\otimes\Z_l
\end{subeqn}
with $\tau$ acting on $N$ by
\begin{subeqn}\label{eqn5412}
\tau\colon (a,b)\mapsto (xa+x\phibar(b)-\phibar(xb),xb)
\end{subeqn}
where
$\phibar\colon\Lambda_{l,r}\otimes\Z_l\to V^1/V^2=\Lambda_{l,r,s}\otimes\Q_l$
is $\phitilde$ composed with the projection $V^1\to V^1/V^2$;
we have $\phibar(x^i)=x^iy_2$ for $0\leq i\leq l^r-2$.

Note that $\Phi_l/\Phi_l^2\cong N/(\tau-1)N$. Hence $\Phi_l/\Phi_l^2$ is cyclic
if and only if the endomorphism $\tau-1$ of the $\F_l$-vector space
$N\otimes\F_l$ has corank $1$. Now $N\otimes\F_l$ is the direct sum of
$\F_l[\eps]/(\eps^{l^{r+s}-l^r})$ and $\F_l[\eps]/(\eps^{l^r-1})$,
with $\eps=x-1$. The matrix of $\tau-1$ with respect to the direct sum of
the bases $(1,\eps,\ldots,\eps^{l^{r+s}-l^r-1})$ and
$(1,\eps,\ldots,\eps^{l^r-2})$ is of the form
\begin{subeqn}\label{eqn5413}
\renewcommand{\baselinestretch}{1}
\left(
\begin{array}{cccc|cccc}
0& & & & & & & \\
1&0& & & &A& & \\
 &\ddots&\ddots& & & & & \\
 & &1&0& & & & \\
\hline
 & & & &0& & & \\
 & & & &1&0& & \\
 & & & & &\ddots&\ddots& \\
 & & & & & &1&0
\end{array}
\right)
\end{subeqn}
It follows that $\tau-1$ has corank $1$ if and only if the upper right
coefficient of $A$ is not zero, or, equivalently, if and only if there
exists $b$ in $\Lambda_{l,r}\otimes\Z_l$ such that $x\phibar(b)-\phibar(xb)$
is a unit in $\Lambda_{l,r,s}\otimes\Z_l$. A computation shows that
$x\phibar(x^{l^r-2})-\phibar(x^{l^r-1})=g_r(x)y_2$. Now recall that we are
free to choose $y=(y_1,y_2)$ in $V^1$ and $z$ in
$(\Lambda_{l,r}\otimes\Q_l)^*$ as long as $(y,z)$ satisfies (\ref{eqn547}).
Note that $g_r(x)$ and $g_r'(x)$ are units in $\Lambda_{l,r,s}\otimes\Q_l$
and $\Lambda_{l,r}\otimes\Q_l$, respectively. Hence we can choose
$y_2=g_r(x)^{-1}$ and $z=x^{-1}g_r'(x)^{-1}$.
\end{example}

\begin{example}\label{example55}
Our final example is the analog of Example~\ref{example54} in the case
$\atilde=0$. More precisely, let $l$ be a prime and $r\geq0$ an integer.
Then there exists an abelian variety $A_K$ with $t=\atilde=0$,
$\ttilde=l^r-1$ and $\Phi=\Phi_l\cong\Z/l^{2r}\Z$.

Let $C_{l,r}$ be as in Example~\ref{example52}. The abelian variety $A_K$
can be found in the isogeny class of $C_{l,r}$ in the same way as used in
Example~\ref{example54}. In this case the construction is somewhat easier,
since the filtration on $V$ has only two steps ($V^1=V^2$), so we leave the
details to the reader. Let us just mention that all formulas up to
(\ref{eqn5410}) remain valid (in adapted form), and after (\ref{eqn5410})
one shows that $M/(\tau-1)M$ can be cyclic with the same method as used to
show that $N/(\tau-1)N$ can be cyclic.
\end{example}

\section{Classification of the $\Phi_{(p)}$.} \label{section6}
The aim of this section is to prove the following theorem.
\begin{theorem}\label{thm61}
Let $D$ be a strictly henselian discrete valuation ring of residue
characteristic $p\geq0$. Let $G$ be a finite commutative group of order
not divisible by $p$. For each prime $l\neq p$, let
$m_l:=(m_{l,1},m_{l,2},\ldots)$ be the partition corresponding to the
$l$-part $G_l$ of $G$ (i.e., $G_l\cong\oplus_{i\geq1}\Z/l^{m_{l,i}}\Z$ and
$m_{l,1}\geq m_{l,2}\cdots$). Let $d$, $t$, $a$ and $u$ be non-negative
integers such that $d=t+a+u$. Then there exists an abelian variety over
the field of fractions of $D$, of dimension $d$, toric rank $t$, abelian
rank $a$ and unipotent rank $u$ which has $\Phi_{(p)}\cong G$, if and
only if
\begin{subeqn}\label{eqn611}
u \geq \sum_{l\neq p}\sum_{i\geq t+1}\left(
\frac{l^{\lfloor m_{l,i}/2\rfloor}+l^{\lceil m_{l,i}/2\rceil}}{2}-1\right)
\end{subeqn}
where for any real number $x$, $\lfloor x\rfloor$ and $\lceil x\rceil$ denote
the largest (resp. smallest) integer $\leq x$ (resp. $\geq x$).
\end{theorem}
\begin{proof}
We will start by showing that if $A_K$ is as indicated in the theorem, then
(\ref{eqn611}) holds. Let $l\neq p$ be a prime. Let $f_l$ be the map from
the set of partitions to $\R$ defined by
\begin{subeqn}\label{eqn612}
f_l(m) = \sum_{i\geq1} \left(\frac{l^{\lfloor m_i/2\rfloor}+
l^{\lceil m_i/2\rceil}}{2}-1\right)
\end{subeqn}
Then $f_l$ is strictly increasing for the partial ordering in which $a\geq b$
if and only if $a_i\geq b_i$ for all~$i$. One easily sees that $f_l$ is
increasing for the lexicographical ordering on the set of partitions of a
fixed number, but we won't use that.
Consider the filtration
\begin{subeqn}\label{eqn613}
\Phi_l \supset \Phi_l^1 \supset \Phi_l^3 \supset 0
\end{subeqn}
induced by (\ref{eqn21}). Theorem~\ref{thm33} shows that
\begin{subeqn}\label{eqn614}
2(t_l-t + a_l-a) \geq \delta_l(\Phi_l/\Phi_l^1) + \delta_l(\Phi_l^1/\Phi_l^3)
\end{subeqn}
Let $n_l$ be the invariant of $\Phi_l/\Phi_l^3$. Lemma~\ref{lemma410}
shows that for all $i\geq1$ we have $n_{l,i}\geq m_{l,i+t}$, or, in the
terminology of the proof of that lemma, that $n_l\geq d^t(m_l)$ in the partial
ordering. Lemma~\ref{lemma411} says that
\begin{subeqn}\label{eqn615}
\delta_l(\Phi_l/\Phi_l^1)+\delta_l(\Phi_l^1/\Phi_l^3)\geq 2f_l(n_l)
\end{subeqn}
It follows that
\begin{subeqn}\label{eqn616}
2(t_l-t+a_l-a) \geq 2f_l(n_l)
\end{subeqn}
Summing over all $l\neq p$ and dividing by $2$ gives (\ref{eqn611}).

It remains to show that all groups $G$ satisfying (\ref{eqn611}) can
occur as the $\Phi_{(p)}$ of an abelian variety $A_K$ over the field of
fractions $K$ of $D$ of dimension $d$, toric rank $t$, abelian rank $a$ and
unipotent rank $u$. It is sufficient to show that all groups $G$ satisfying
\begin{subeqn}\label{eqn617}
u = \lceil\left(\sum_{l\neq p}\sum_{i\geq t+1}\left(
\frac{l^{\lfloor m_{l,i}/2\rfloor}+l^{\lceil m_{l,i}/2\rceil}}{2}-1\right)
\right)\rceil
\end{subeqn}
occur in such a way, since one can replace $A_K$ by the product of $A_K$ with
an abelian variety $B_K$ which has unipotent reduction and trivial group
of connected components.

Let us first suppose that $K=\C((q))$. Let $d$, $t$, $a$, $u$ and $G$ be as
in the theorem, and suppose that they satisfy (\ref{eqn617}). We have
$G\cong \oplus_{i\geq1}\Z/n_i\Z$ with $n_i\geq1$ and $n_{i+1}|n_i$ for all
$i$. Let $B_K$ be of the type described in Example~\ref{example51}:
it has dimension $t$, completely toric reduction and group of connected
components $\Z/n_1\Z\oplus\cdots\oplus\Z/n_t\Z$. The abelian variety $A_K$
we are constructing will be of the form
\begin{subeqn}\label{eqn618}
A_K=B_K\times\prod_{l\neq p}C_{K,l}, \qquad{\rm with}\quad
\dim(C_{K,l})=\lceil\left(\sum_{i\geq t+1}\left(
\frac{l^{\lfloor m_{l,i}/2\rfloor}+l^{\lceil m_{l,i}/2\rceil}}{2}-1\right)
\right)\rceil
\end{subeqn}
and such that all $C_{K,l}$ have unipotent reduction. Note that in fact
such an $A_K$ has unipotent rank $u$, since for $l\neq2$ the function $f_l$
defined above has integer values. For $l\neq2$ we define
\begin{subeqn}\label{eqn619}
C_{K,l} = \prod_{i>t} C_{K,l,i}
\end{subeqn}
where $C_{K,l,i}$ is the abelian variety constructed in Example~\ref{example54}
with $r=(m_{l,i}-1)/2$ and $s=1$ if $m_{l,i}\neq1$ is odd, where $C_{K,l,i}$
is the abelian variety constructed in Example~\ref{example53} with $i=1$ if
$m_{l,i}=1$, and $C_{K,l,i}$ is the abelian variety constructed in
Example~\ref{example55} with $r=m_{l,i}/2$ if $m_{l,i}$ is even.
Note that the group of connected components of
the reduction of $C_{K,l,i}$ is cyclic of order~$l^{m_{l,i}}$. For $l=2$ the
construction of $C_{K,l}$ is a bit different. Let $r\geq0$ be maximal such
that $m_{2,i}>1$ for all $i\leq r$. Then $C_{K,2}$ will be of the form
\begin{subeqn}\label{eqn6110}
C_{K,2} = D_{K,2} \times \prod_{t<i\leq r}C_{K,2,i}
\end{subeqn}
where $C_{K,2,i}$ is defined as $C_{K,l,i}$ but with $l$ replaced by $2$, and
where $D_{K,2}$ is as follows. Let $v$ be the number of $i>t$ such that
$m_{2,i}=1$. If $v$ is even we let $D_{K,2}$ be the product of $v/2$
elliptic curves which have unipotent reduction and group of connected
components isomorphic to $\Z/2\Z\times\Z/2\Z$. If $v$ is odd we let
$D_{K,2}$ be the product of $(v-1)/2$ elliptic curves with unipotent reduction
and group of connected components isomorphic to $\Z/2\Z\times\Z/2\Z$ and one
elliptic curve with unipotent reduction and group of connected components
cyclic of order~$2$. One verifies easily that $A_K$ has all the desired
properties.

To finish the proof of the theorem, we have to show that similar examples
exist over any strictly henselian discrete valuation ring $D$ with residue
characteristic $p$. Since our examples are products of the examples of
\S\ref{section5}, it suffices to show that the examples in \S\ref{section5}
exist over $D$. Since we do not suppose $D$ complete, we cannot use a
Tate curve ``$\Gm/q^\Z$'' with $q$ a uniformizer of~$D$. Instead we can
use any elliptic curve $E$ over $K$ which has toric reduction and trivial
group of connected components. Then $I$ acts on the Tate module
$\rT_l(E(\Ksep))$ through its quotient $\Z_l(1)$ and for a suitable choice
of a $\Z_l$-basis of $\rT_l(E(\Ksep))$, an element of $I$ with image $a$
in $\Z_l(1)$ acts as $({1\atop0}{a\atop1})$. It follows that
Examples~\ref{example51} and \ref{example52} with $E_1$ replaced by $E$ still
work. To make Example~\ref{example53} work over $D$, it is enough to show
that for all $l\neq p$ and $r>0$ such that $l^r>2$, there exists an abelian
scheme over $D$ of relative dimension $l^{r-1}(l-1)/2$ and with an action by
$\Z[x]/(f_{l,r})$. Once one has these abelian schemes, the constructions
of \S\ref{section5} can be carried out over~$D$. The fact that such abelian
schemes exist is a consequene of the theory of abelian varieties of
``CM-type''. Fix an $l$ and $r$ as above. The moduli scheme over $\Z[1/l]$ of
abelian schemes of the  desired type, with a suitable polarization and
$l$-power level structure, is finite etale and not empty. Another way to
prove the desired existence is to consider isogeny factors over
$\Q(\zeta_{l^r})$ of the jacobian of the Fermat curve of degree $l^r$.
\end{proof}

\section{Further remarks and questions.} \label{section7}
Although Theorem~\ref{thm61} gives a complete classification of the
prime-to-$p$ parts of the groups of connected components of special fibres of
N\'eron models with some fixed invariants, there are still questions left.
For example, it is clear
that the groups of connected components $\Phi$ have functorial additional
structure coming from the fact that the category of abelian varieties has
an involution: every abelian variety has its dual. More precisely, suppose
that $\Phi$ comes from the abelian variety $A_K$. Let $A'_K$ be the
dual of $A_K$ and denote its group of connected components by $\Phi'$.
Then there are several pairings with values in $\Q/\Z$, conjecturally perfect
and the same, between $\Phi$ and $\Phi'$; see
\cite[\S\S1.2, 1.3, 11.2]{Grothendieck1}, \cite{Milne1}, \cite[\S3]{Dino1},
\cite{Moret-Bailly1} and
\cite[Prop.~3.3]{Oort1}. Let us note by the way that the last reference is
clearly wrong since it says that the pairing has values in $(\Q/\Z)(1)$;
the mistake in the proof is that the direct sum decomposition in the unique
displayed formula in it is not unique. Anyway, for each of the remaining
pairings we get a filtration
\begin{eqn}\label{eqn71}
\Phi_{(p)} = {\Phi'}_{(p)}^{4,\perp} \supset {\Phi'}_{(p)}^{3,\perp} \supset
{\Phi'}_{(p)}^{2,\perp} \supset {\Phi'}_{(p)}^{1,\perp} \supset
{\Phi'}_{(p)}^{0,\perp} = 0
\end{eqn}
It would be interesting to know the common refinement of this filtration
with (\ref{eqn21}). Also, it would be of interest to prove that the various
pairings are the same up to a determined sign. Some relations between
the two filtrations (\ref{eqn21}) and (\ref{eqn71}) on the $l$-part for
$l\neq p$ can be found in \cite[Thm.~3.21]{Dino1}, under the hypothesis that
$A_K$ has a polarization of degree prime to~$l$.

Let us consider the functor from the category of abelian varieties over $K$
to the category of finite abelian groups which associates to each abelian
variety the group of connected components of the special fibre of its
N\'eron model. A rather vague question one can ask is through what categories
of abelian groups endowed with some extra structure this functor factors.
We have seen for example that there is a filtration of four steps on the
prime-to-$p$ part, but as we have just remarked that is certainly not all
there is.

Lorenzini has shown \cite[Thm.~3.22]{Dino1}, under the hypothesis that there
is a polarization of degree prime to $l$, that ${\Phi'}_l^{2,\perp}$ is the
prime-to-$p$ part of the kernel of the map from $\Phi$ to the group of
connected components of $A_L$, where $K\to L$ is any extension over which
$A_K$ has semi-stable reduction. It would be interesting to generalize this.
Even in the case in which $A_K$ acquires semi-stable reduction after a
tamely ramified extension $K\to L$, when the theory of \cite{Edixhoven1}
applies, I have not been able to give a description of the filtration
(\ref{eqn21}) in terms of the special fibre of the N\'eron model of $A_L$
with its action of $\Gal(K/L)$.

The $p$-part of $\Phi$ remains difficult. For example, one expects a bound
for its order in terms of the dimension of $A_K$ if the toric part of the
reduction is zero, but even in the case of potentially good reduction I don't
know of any such bound (of course, if $A_K$ is the jacobian of a curve with a
rational point, the usual bound, i.e., the bound we have when $k$ is of
characteristic zero, holds, since one can apply Winters's theorem
\cite{Winters1}).
In a forthcoming article \cite{EdiLiuLor} one can
find a generalization of a result of McCallum (unpublished) which says that
in the case of potentially good reduction the $p$-part is annihilated by the
degree of any extension after which one obtains semi-stable reduction, but
not in general by the exponent of the Galois group of such an extension.
Work in progress by Bosch and Xarles, using a rigid analytic uniformization of
N\'eron models, seems to imply that there is a four-step functorial filtration
on the whole of~$\Phi$, for which three of the four successive quotients can
be described in terms of the ineria group acting on the character group of the
toric part of the semi-stable reduction. The remaining successive quotient
comes from an abelian variety, obtained by Raynaud's extension, which has
potentially good reduction. This part is still a mystery.

\vspace{\baselineskip}
\noindent
Institut Math\'ematique \\
Universit\'e de Rennes 1 \\
Campus de Beaulieu \\
35042 Rennes cedex \\
France.


\begin{thebibliography}{99}

\bibitem{BLR} S.\ Bosch, W.\ L\"utkebohmert and M.\ Raynaud.
{\sl N\'eron models.} Ergebnisse der Mathematik und ihrer Grenzgebiete,
3.~ Folge, Band~21, Springer-Verlag (1990).

\bibitem {DeligneRapoport} P.~Deligne and M.~Rapoport. {\sl Les sch\'emas
de modules de courbes elliptiques.} Modular Functions of One Variable II,
143--316. Lecture Notes in Mathematics 349, Springer-Verlag (1975).

\bibitem {Edixhoven1} B.\ Edixhoven. {\sl N\'eron models and tame
ramification.} Compositio Mathematica~81 (1992), 291--306.

\bibitem {EdiLiuLor} B.\ Edixhoven, Q.\ Liu and D.\ Lorenzini.
{\sl The $p$-part of the group of components.} To be submitted for
publication.

\bibitem {Grothendieck1} A.\ Grothendieck. {\sl Mod\`eles de N\'eron et
monodromie.} S\'eminaire de G\'eom\'etrie~7, Expos\'e~IX, Lecture
Notes in Mathematics 288, Springer-Verlag (1973).

\bibitem {LenstraOort} H.W.\ Lenstra and F.\ Oort. {\sl Abelian varieties
having purely additive reduction.} J.~Pure Appl.\ Algebra~36 (1985),
281--298.

\bibitem {Dino1} D.\ Lorenzini. {\sl On the group of components of a N\'eron
model.} J.~reine angew.~Math.~445 (1993), 109--160.

\bibitem {MacDonald1} I.G.\ MacDonald. {\sl Symmetric functions and Hall
polynomials.} Oxford, Clarendon Press (1979).

\bibitem {Milne1} J.S.\ Milne. {\sl Arithmetic duality theorems.} Perspectives
in Mathematics~1, Academic Press (1986).

\bibitem {Moret-Bailly1} L.\ Moret-Bailly. {\sl Pinceaux de vari\'et\'es
ab\'eliennes.} Ast\'erisque~129. Soci\'et\'e Math\'ematique de France (1985).

\bibitem {Mumford} D.~Mumford. {\sl An analytic construction of degenerating
abelian varieties over complete rings.} Compositio Mathematica~24, Fasc.~3
(1972), 239--272.

\bibitem {Oort1} F.\ Oort. {\sl Good and stable reduction of abelian
varieties.} Manuscripta Math.~11 (1974), 171--197.

\bibitem {Tate1} J.\ Tate. {\sl Algorithm for determining the type of a
singular fibre in an elliptic pencil.} Modular Functions of One
Variable IV, 33--52. Lecture Notes in Mathematics 476, Springer-Verlag (1975).

\bibitem {Winters1} G.\ Winters. {\sl On the existence of certain families of
curves.} American Journal of Mathematics~96 (1974), 215--228.


\end{thebibliography}
\end{document}